\documentclass[a4paper, 12pt]{amsart}
\usepackage{amsmath,mathrsfs,amsthm}
\usepackage{txfonts,rotating,tikz}
\usetikzlibrary{intersections,shapes.arrows,calc}

\usepackage[all,cmtip]{xy}
\usepackage{amssymb}
\usepackage{amsfonts}
\usepackage{amscd}
\usepackage{color} 
\usepackage[utf8]{inputenc}
\usepackage[headings]{fullpage}
\usepackage{epic}
\usepackage{mathtools}
\usepackage{tikz-cd}
\usetikzlibrary{chains}
\usepackage[thinlines]{easytable}
\usepackage{array}
\usepackage{ctable}
\usepackage{pbox}
\usepackage[usestackEOL]{stackengine}
\usepackage[sort, numbers]{natbib}
\usepackage{hyperref}
\usepackage[nameinlink,noabbrev]{cleveref}
\usepackage{booktabs,tabularx,array}
\usepackage{float}
\tikzcdset{scale cd/.style={every label/.append style={scale=#1},
 cells={nodes={scale=#1}}}}

\newcolumntype{?}{!{\vrule width 1pt}}

\newtheorem{introconjecture}{Conjecture}

\newtheorem{introtheorem}[introconjecture]{Theorem}

\newtheorem*{theorem*}{Theorem}

\newtheorem{theorem}{Theorem}[section]

\newtheorem{lemma}[theorem]{Lemma}
\newtheorem{proposition}[theorem]{Proposition}
\newtheorem{corollary}[theorem]{Corollary}

\newtheorem{remark}[theorem]{Remark}

\theoremstyle{definition}
\newtheorem{definition}[theorem]{Definition}

\newtheorem{example}[theorem]{Example}

\newcommand{\be}{\begin{equation}}
	\newcommand{\ee}{\end{equation}}
\newcommand{\bt}{\begin{theorem}}
	\newcommand{\et}{\end{theorem}}
\newcommand{\bd}{\begin{definition}}
	\newcommand{\ed}{\end{definition}}
\newcommand{\bp}{\begin{proposition}}
	\newcommand{\ep}{\end{proposition}}
\newcommand{\bl}{\begin{lemma}}
	\newcommand{\el}{\end{lemma}}
\newcommand{\bco}{\begin{corollary}}
	\newcommand{\eco}{\end{corollary}}
\newcommand{\br}{\begin{remark}}
	\newcommand{\er}{\end{remark}}
\newcommand{\bex}{\begin{example}}
	\newcommand{\eex}{\end{example}}
\newcommand{\ben}{\begin{enumerate}}
	\newcommand{\een}{\end{enumerate}}
\newcommand{\bc}{\begin{cases}}
	\newcommand{\ec}{\end{cases}}

\newcommand{\bpf}{\begin{proof}}
	\newcommand{\epf}{\end{proof}}

\newcommand{\bma}{\begin{bmatrix}}
	\newcommand{\ema}{\end{bmatrix}}

\newcommand{\C}{\mathbb C}
\newcommand{\Z}{\mathbb Z}
\newcommand{\g}{\mathfrak g}

\newcommand{\fh}{\mathfrak{h}}

\newcommand{\la}{\langle}
\newcommand{\ra}{\rangle}

\DeclareMathOperator{\Gr}{\mathrm{Gr}}

\definecolor{BHYblue}{RGB}{37,76,116}
\definecolor{BHYgray}{RGB}{244,246,248}
\newcommand{\beqn}{\begin{equation}}
	\newcommand{\eeqn}{\end{equation}}

\newcommand{\arxiv}[1]{\href{https://arxiv.org/abs/#1}{\texttt{arXiv:#1}}}

\begin{document}
\title{Nonreduced Finkelberg-Mirkovi\'c Schubert schemes}
\author{Marc Besson}
	\address[M.Besson]{Shanghai Institute for Mathematics and Interdisciplinary Science, Shanghai, China}
\email{bessonm@simis.cn}
	\author{Jiuzu Hong}
	\address[J.Hong]{Department of Mathematics, University of North Carolina at Chapel Hill, Chapel Hill, NC 27599-3250, U.S.A.}
	\email{jiuzu@email.unc.edu}
	\author{Huanhuan Yu}
	\address[H.Yu]{Department of Mathematics, Jiangxi University of Finance and Economics, Nanchang, 330032, China}
	\email{yuhuanhuan@jxufe.edu.cn}

\begin{abstract}

We develop an effective numerical criterion ($\theta$-minimality) to detect the nonreducedness of Finkelberg-Mirkovi\'c (FM) Schubert schemes  via tangent-space comparisons at arbitrary torus-fixed points.  We completely classify the dominant coweights admitting $\theta$-minimal pairs.  As an application, we obtain nonreduced FM Schubert schemes in every type other than $A_n$ and $C_n$. In particular, every nontrivial FM Schubert scheme of type $E_8$ is nonreduced. In type $D_n$, we show that some nonreduced FM Schubert schemes and their corresponding spherical Schubert varieties have the same tangent space at the base point. Our methods use the Frenkel-Kac construction of basic representations and Segal-Sugawara operators.
\end{abstract}
\maketitle

\setcounter{tocdepth}{2}
\tableofcontents

\section{Introduction}

Let $G$ be a simple algebraic group over $\C$ with
maximal torus $T$. Let $\Gr_G$ be the affine Grassmannian of $G$.  For a
dominant coweight $\lambda$, the spherical Schubert variety
$\Gr_{\preceq\lambda}$ is the reduced closure of the $L^+G$-orbit
through the $T$-fixed point $e_\lambda\in\Gr_G$. 

When $G$ is simply-connected, Finkelberg and Mirkovi\'c \cite{FM} proposed a moduli-theoretic
description of $\Gr_{\preceq\lambda}$, by imposing pole conditions for  associated lattices attached to all irreducible
representations of $G$ (see Section~\ref{sec:fmtest}).  The associated 
\emph{Finkelberg-Mirkovi\'c (FM) Schubert scheme}
$\Gr^{\rm FM}_{\preceq\lambda}$ is set-theoretically supported on
$\Gr_{\preceq\lambda}$, but a priori it may not be reduced.  More precisely, there is a
closed immersion  $\Gr_{\preceq\lambda}\hookrightarrow\Gr^{\rm FM}_{\preceq\lambda}$ 
identifying $\Gr_{\preceq\lambda}$ with the reduced subscheme of
$\Gr^{\rm FM}_{\preceq\lambda}$.  The \emph{reducedness conjecture}
(cf.\,\cite{KMW,KWWY}) asserts that this immersion is an
isomorphism for every $\lambda$.  It was proved in type $A$ by
Kamnitzer, Muthiah, Weekes and Yacobi \cite{KMWY}.  The reducedness question is related to
the quantization of affine Grassmannian slices by truncated shifted Yangians
\cite{KWWY}, as well as the construction of integral models of Shimura varieties \cite{KPZ}.
For a general simple algebraic group $G$, following \cite{ZhuPCMI,KPZ} we use a
$z$-extension of $G$ to formulate the FM Schubert scheme uniformly in every
connected component of $\Gr_G$. In Proposition~\ref{prop:zfm-independence}, we show that
the resulting closed subscheme is independent of the choice of $z$-extension.

In our previous work \cite{BHY}, we approached the reducedness conjecture by
comparing the Zariski tangent spaces of $\Gr_{\preceq\lambda}$ and
$\Gr^{\rm FM}_{\preceq\lambda}$.  In particular, we determined the
tangent space of every FM Schubert scheme at the base point
$e\in\Gr_G$, which led to unexpected counterexamples to the reducedness conjecture
in types $E_6$, $E_7$ and $E_8$.  The simplest of our counterexamples is the
quasi-minuscule Schubert variety in type $E_8$, for which
\[
 \dim T_e\Gr_{\preceq\theta^\vee}=\dim G,
 \qquad
 \dim T_e\Gr^{\rm FM}_{\preceq\theta^\vee}=2\dim G.
\]
Here $\theta^\vee$ is the coroot of the highest root $\theta$.
All counterexamples in \cite{BHY} were detected by comparing tangent spaces at $e$. 
\vspace{2mm}

In this work, we study the tangent spaces of $\Gr_{\preceq \lambda}$ and $\Gr_{\preceq \lambda}^{\rm FM}$ at arbitrary $T$-fixed points
$e_\mu$.  This flexibility leads to a much better understanding of the difference of these tangent spaces. In particular, we obtain counterexamples to the reducedness conjecture for every simple group of type $B_n (n\geq3), D_n (n\geq 4), E_6, E_7, E_8, F_4, G_2$.  The method used here is different from that of \cite{BHY}. In addition to Polo's criterion, which was already used  in \cite{BHY}, in the simply-laced case we employ the Frenkel-Kac construction of basic representations together with Segal-Sugawara operators to show that certain tangent vectors of $T_{e_\mu}\Gr^{\rm FM}_{\preceq \lambda}$ do not belong to $T_{e_\mu}\Gr_{\preceq \lambda}$; for non simply-laced types, we use diagram folding.

\vspace{1mm}
For the statements of our main theorems, we now assume that $G$ is adjoint.
Let $Q,P,Q^\vee,P^\vee$ denote the root, weight, coroot, and coweight
lattices, respectively. Thus, $X^*(T)=Q$ and $X_*(T)=P^\vee$.
We write $X_*(T)^+$ for the dominant coweights and
$\langle\ ,\ \rangle$ for the natural pairing. Let $\theta$ be the
highest root, and let $\alpha_i$, $\alpha_i^\vee$, $\omega_i$, and
$\omega_i^\vee$ denote the simple roots, simple coroots, fundamental
weights, and fundamental coweights, respectively. For
$\mu,\lambda\in X_*(T)^+$, we write $\mu\preceq\lambda$ if
\[
\lambda-\mu=\sum_i a_i\alpha_i^\vee,
\qquad a_i=\langle\lambda-\mu,\omega_i\rangle\in\mathbb Z_{\geq0},
\]
and $\mu\prec\lambda$ if moreover $\mu\neq\lambda$.

For $\mu\prec\lambda$, we say that the pair $(\lambda,\mu)$ is $\theta$-\emph{minimal} if
\begin{equation}\label{eq:intro-admissible}
1\leq\langle\lambda-\mu,\theta\rangle
\leq\min_i a_i.
\end{equation}

The following theorem is one of the main results of this paper. 
\begin{introtheorem}[Theorem \ref{thm:main1}]\label{thm:intro-main}
Let $\lambda\in X_*(T)^+$ be a dominant coweight. Suppose that there exists
$\mu\in X_*(T)^+$ such that $(\lambda,\mu)$ is $\theta$-minimal.  Then
\[
 T_{e_\mu}\Gr_{\preceq\lambda}
 \subsetneq
 T_{e_\mu}\Gr^{\rm FM}_{\preceq\lambda}.
\]
As a consequence, the FM Schubert scheme $\Gr^{\rm FM}_{\preceq\lambda}$ is nonreduced.
\end{introtheorem}

There are no $\theta$-minimal pairs for $G$ of type $A$, which is consistent with
\cite{KMWY}. There are also no $\theta$-minimal pairs in types $C$ and $B_2$.
In types $B_n$ ($n\geq3$), $D$, $E$, $F_4$ and $G_2$, the highest root $\theta$ is a fundamental
weight $\omega_s$, where $s$ is the corresponding node. So $(\lambda,\mu)$ is $\theta$-minimal if and only
if $a_s=\min_i a_i\geq1$ (Lemma~\ref{lem:criterion}).
We completely classify the dominant coweights admitting a
$\theta$-minimal pair in every simple type.
\begin{samepage}
\begin{introtheorem}[Theorem \ref{thm:classification-all}]\label{thm:intro-classification}
Write $\lambda=\sum_i b_i\omega_i^\vee\in X_*(T)^+$. There exists $\mu\in X_*(T)^+$ such that $\mu\prec\lambda$ and $(\lambda,\mu)$ is $\theta$-minimal exactly in the following cases:
\begin{center}
\small
\renewcommand{\arraystretch}{1.18}
\begin{tabularx}{0.98\textwidth}{@{}l@{\hspace{2em}}X@{}}
\toprule
Type & Condition on $\lambda$ \\
\midrule
$A_n$, $C_n$, $B_2$ & none \\
$B_n$, $n\ge3$ & $b_1\ge1$ and $\sum_{i=3}^n b_i\ge1$ \\
$D_n$, $n\ge5$ & $b_1\ge1$ and $2\sum_{i=3}^{n-2}b_i+b_{n-1}+b_n\ge2$ \\
$D_4$ & $b_1,b_3,b_4\ge1$ \\
$E_6$ & $5b_1+4b_3+3b_4+2b_5+b_6\ge3$ and $b_1+2b_3+3b_4+4b_5+5b_6\ge3$ \\
$E_7$ & $b_2+b_7+2(b_3+b_4+b_5+b_6)\ge2$ \\
$E_8$ & $\lambda\neq0$ \\
$F_4$ & $b_2+b_3+b_4\ge1$ \\
$G_2$ & $b_1\ge1$ \\
\bottomrule
\end{tabularx}
\end{center}
\end{introtheorem}
\end{samepage}
In particular, every nontrivial FM Schubert scheme of type $E_8$ is nonreduced (Corollary~\ref{cor:E8-all}).  The reducedness theorem of \cite{KMWY} is proved for $\mathrm{SL}_n$ and does not address FM Schubert schemes for $\mathrm{PGL}_n$ in non-neutral components. The reducedness conjecture in type $C_n$ also remains open in general.
\vspace{1mm}

Spherical Schubert varieties associated to sums of minuscule coweights arise naturally in the geometry of local models of Shimura varieties.
In \cite[Theorem~4.2.3]{KPZ}, Kisin, Pappas and Zhou considered
$\lambda=r\omega_1^\vee$ in types $B_n$ and $D_n$, and
$\lambda=s\omega_{n-1}^\vee+t\omega_n^\vee$ in type $D_n$.
For these coweights, they showed that the tangent spaces of
$\Gr_{\preceq\lambda}$ and $\Gr^{\rm FM}_{\preceq\lambda}$
agree in all real root directions at every $T$-fixed point.
These families considered in \cite{KPZ} lie in the complementary range of the classification in Theorem \ref{thm:intro-classification}.
\vspace{1mm}

Another new phenomenon which appears in type $D_n$ is the following. Proposition \ref{prop:no-basepoint-D} shows that in type $D_n$, there is no $\theta$-minimal pair $(\lambda, \mu)$ with $\mu=0$. However, there may still exist $\theta$-minimal pairs $(\lambda, \mu)$ for $\mu \neq 0$. In fact, we have:
\begin{samepage}
\begin{introtheorem}[Theorem \ref{lem:nonbase}]
For every type \(D_n\), there exist dominant coweights $\lambda\in X_*(T)^+\cap Q^\vee$ such that
\[
T_e\Gr_{\preceq\lambda}
=
T_e\Gr^{\rm FM}_{\preceq\lambda},
\]
while $(\lambda, \theta^\vee)$ is a $\theta$-minimal pair. In particular,  \(\Gr^{\rm FM}_{\preceq\lambda}\) is nonreduced.
\end{introtheorem}
\end{samepage}
 This theorem gives a negative answer to a question raised by Alex Weekes: equality of tangent spaces at $e$ does not  imply that the FM Schubert scheme is reduced.  Thus, in type \(D_n\), passing to non-base fixed points is essential for detecting the nonreducedness of these FM Schubert schemes. The proof of this theorem uses the results of \cite{BHY}
and the Frenkel-Kac construction of basic representation. 

\vspace{2mm}
We now explain the idea of the proof of Theorem \ref{thm:intro-main}. Let $x_\theta\in\g$ be a root vector for the highest
root $\theta$, put $q=\langle\lambda-\mu,\theta\rangle$, and consider the
tangent vector $\xi_{\mu,q}\in T_{e_\mu}\Gr_G$ at $e_\mu$ corresponding to the Lie algebra element $x_\theta(-q):=x_\theta\otimes t^{-q}$. If $(\lambda, \mu)$ is $\theta$-minimal,
we show that $\xi_{\mu,q}$ is tangent to $\Gr^{\rm FM}_{\preceq\lambda}$ but not to
$\Gr_{\preceq\lambda}$. 

For simply-laced $G$, we identify coweights with weights via the normalized
invariant form in this introduction, and write $D(1,\lambda)$ for the
corresponding level-one Demazure module.

A key feature of our approach is the use of Segal-Sugawara operators to detect tangent vectors that do not belong to $T_{e_\mu}\Gr_{\preceq \lambda}$.  By Polo's criterion, this amounts to
showing that a certain vector of the level one representation does not belong to the Demazure module $D(1,\lambda)$.  To determine whether a given vector lies in a Demazure module is difficult in
general. We approach this by using the Segal-Sugawara operators $L_k$. The key observation is that $L_k$ annihilates the whole Demazure module when $k$ is large relative to $\lambda$.
\begin{introtheorem}[Theorem \ref{thm:vir-bound}]
Let $G$ be simply-laced. For any dominant coweight $\lambda\in X_*(T)^+$ and  every integer  $k\geq \max ( \langle\lambda,\theta\rangle -1, 1)$,  we have 
\begin{equation}
L_k \cdot D(1,\lambda)=0, 
\end{equation}
where $L_k$ is the $k$-th Segal-Sugawara operator. 
\end{introtheorem}

For $q=\langle\lambda-\mu,\theta\rangle\geq1$, we describe a positive
Segal--Sugawara operator that annihilates $D(1,\lambda)$ but acts
nontrivially on $x_\theta(-q)v_{-\mu}$. This proves
$x_\theta(-q)v_{-\mu}\notin D(1,\lambda)$ (Proposition~\ref{prop:nonincl}), where $v_{-\mu}$ is the extremal
weight vector corresponding to $e_\mu$. By Polo's criterion, the
tangent vector $\xi_{\mu,q}$ therefore does not belong to
$T_{e_\mu}\Gr_{\preceq\lambda}$. On the other hand, in Proposition~\ref{prop:fmtest} we show that the defining pole conditions of FM Schubert scheme give
\[
\xi_{\mu,q}\in T_{e_\mu}\Gr^{\rm FM}_{\preceq\lambda}
\quad\Longleftrightarrow\quad
q\leq\min_i\langle\lambda-\mu,\omega_i\rangle.
\]
The nonsimply-laced case follows by realizing the adjoint group as the
fixed subgroup of a pinned automorphism of a simply-laced adjoint group
and applying the tangent-space inclusion induced by the embedding of
affine Grassmannians, using the idea from \cite{HongMV}. Both steps are included in the proof of
Theorem~\ref{thm:main1}.

 Thus, for $q=\langle\lambda-\mu,\theta\rangle$, the $\theta$-minimality condition ensures that
$\xi_{\mu,q}\in T_{e_\mu}\Gr^{\rm FM}_{\preceq\lambda}$, but
$\xi_{\mu,q}\notin T_{e_\mu}\Gr_{\preceq\lambda}$. Hence every $\theta$-minimal pair yields a tangent vector detecting the nonreducedness of the FM Schubert scheme, proving the main theorem.

 \vspace{0.2cm}

{\bf Acknowledgements}:   We would like to thank Alex Weekes for raising a question in email correspondence that led to Theorem \ref{lem:nonbase}. A significant part of this work was carried out while M. Besson was a Shuimu postdoctoral researcher at YMSC, Tsinghua. J.\,Hong gratefully acknowledges support from Simons
Travel Grant MPS-TSM-00007468.  H.\,Yu was supported by the National Natural Science Foundation of China (Grant No. 12601054).

\section{Demazure modules and Sugawara operators}

In this section, we recall the Frenkel-Kac realization of level-one representations
and the Segal-Sugawara operators.
We then establish a vanishing theorem for these operators on Demazure modules. This gives us an effective criterion to show that many candidate vectors do not belong to affine Demazure modules; in Section \ref{sect:main_thm} we use this to study tangent spaces of spherical Schubert varieties.
\label{sec:FK}

\subsection{Frenkel-Kac construction of  basic representations}\label{subsect:FK}

We recall the Frenkel-Kac construction, following
\cite[Chapter~14]{Kac} closely. Let $\hat{\mathfrak g}$ be the
untwisted affine Kac--Moody algebra associated with a
simply-laced simple Lie algebra $\mathfrak g$, and let
$\mathfrak h\subset \mathfrak g$ be a Cartan subalgebra. Let
$\Phi\subset\mathfrak h^*$ be the root system, with positive roots
$\Phi^+$, simple roots $\alpha_1,\ldots,\alpha_\ell$, and highest root
$\theta$. We use $(\,,\mid\,,)$ for the invariant bilinear form
on $\mathfrak g$ and its induced form on $\mathfrak h^*$, normalized so
that $(\alpha\mid \alpha)=2$ for every root $\alpha$.
Let $Q$ and $P$ denote the root and weight lattices, respectively, and
let $P^+$ be the set of dominant integral weights. Via the normalized
invariant form, we also regard $P$ and $Q$ as lattices in $\mathfrak h$;
in particular, a root $\alpha$ corresponds to its coroot $\alpha^\vee$,
and we use $\alpha$ for this Cartan element as well.
For $\lambda,\mu\in P$, we write $\mu\preceq\lambda$ when
$\lambda-\mu\in\sum_i a_i \alpha_i$ for $a_i \in \mathbb{Z}^{\geq 0}$, and $\mu\prec\lambda$
when in addition $\mu\neq\lambda$.

The affine Kac--Moody algebra is
$\hat{\mathfrak g} = \mathfrak g[t,t^{-1}] \oplus \mathbb CK \oplus \mathbb Cd$,
with Lie bracket
\[
[x(m),y(n)]
=
[x,y](m+n)
+
m\delta_{m,-n}(x\mid y) K,
\qquad
[d,x(n)]=n\,x(n),
\]
where 
$x(n):=x\otimes t^n$,  $K$ is central, and $d$ is the degree operator.

The Frenkel-Kac construction gives a concrete realization of all level one highest weight integrable representations of $\hat{\g}$ in the simply-laced case.  We consider the following vector space:
\begin{equation}
\label{eq:basic-space}
V
=
S\bigl(t^{-1}\mathfrak h[t^{-1}]\bigr)\otimes_{\mathbb{C}} \mathbb C[P],
\end{equation}
where $P\subset \mathfrak h$ and $\mathbb C[P]$ is its group
algebra.  We write $e^\alpha$ for the image of $\alpha$ under the natural inclusion 
\[P \rightarrow \mathbb{C}[P].\]

If $u\in\mathfrak h$ and $n>0$, then $u(-n)$ acts on $V$ by
left multiplication.  For $u \in \mathfrak{h}$ and $n\ge 0$, we define an operator $u(n)$ on $V$ by the formula 
\begin{equation}
\label{eq:heisenberg-action}
u(n)(f\otimes e^\alpha)
=
(u(n)f)\otimes e^\alpha
+
\delta_{n,0}(\alpha\mid u)  (f\otimes e^\alpha).
\end{equation}
More explicitly, if
\[
f=v_1(m_1)\cdots v_k(m_k),
\qquad m_i<0,
\]
then for $n>0$,
\[
u(n)f
=
\sum_{i=1}^k
n\,\delta_{n,-m_i}(u\mid v_i)
\,v_1(m_1)\cdots \widehat{v_i(m_i)}\cdots v_k(m_k).
\]

Choose an asymmetry function $\epsilon:Q\times Q\to \{\pm1\}$ as in \cite[\S7.8]{Kac}, and extend it bimultiplicatively to
$\epsilon:P\times P\to \mathbb C^\times$,
using the same notation. For $\beta\in Q$, define the operator
$c_\beta$ by
\[
c_\beta(f\otimes e^\lambda)
=
\epsilon(\beta,\lambda)f\otimes e^\lambda,
\qquad \lambda\in P.
\]
Then $\epsilon(\beta,0)=1$ for every $\beta\in Q$, and $\epsilon(\alpha,\alpha)=\epsilon(\alpha,-\alpha)=-1$
for every root $\alpha$.

For each root $\alpha$, let $\mathfrak g_\alpha$ denote the
corresponding root space in $\g$, and choose a nonzero root vector
$x_\alpha\in\mathfrak g_\alpha$ as in \cite[\S 7.8]{Kac} satisfying the properties in \cite[(7.8.6)]{Kac}. In particular, 
we have 
\begin{equation}\label{eq:x-bracket}
 (x_\alpha\mid x_{-\alpha})=-1, \quad [x_\alpha,x_{-\alpha}]=-\alpha  ,\end{equation}
where $\alpha$ is identified with the corresponding element of
$\mathfrak h$ via the invariant form.

For a root
$\alpha\in\Phi$, define the vertex operator
\begin{equation}
\label{eq:Gamma-alpha}
x_\alpha(z)
=
\exp\!\left(\sum_{j\ge1}\frac{z^j}{j}\alpha(-j)\right)
\exp\!\left(-\sum_{j\ge1}\frac{z^{-j}}{j}\alpha(j)\right)
e^\alpha z^{\alpha(0)}c_\alpha.
\end{equation}
We write
\[
x_\alpha(z)
=
\sum_{n\in\mathbb Z}\Gamma_\alpha(n)z^{-n-1}.
\]

Let $u_i,u^i$ be dual bases of $\mathfrak h$ under the normalized invariant form $(\,\mid\,)$.  Define
\begin{equation}
\label{eq:D0}
D_0
=
\frac12\sum_{i=1}^{\ell}u_i(0)u^i(0)
+
\sum_{i=1}^{\ell}\sum_{n\ge1}u_i(-n)u^i(n).
\end{equation}

The following theorem, due to Frenkel-Kac \cite{FK}, and Segal \cite{Se}, also see \cite{Do}, realizes all level-one integrable highest weight
representations on $V$ via the operators defined above. 
\begin{theorem}\label{thm:FKS}
For $\hat{\mathfrak{g}}$ an affine algebra of type A, D or E, the map $\hat{\mathfrak{g}} \rightarrow {\rm End}(V)$
\[
K\longmapsto 1,
\qquad
u\otimes t^n\longmapsto u(n)
\quad (u\in\mathfrak h),
\qquad
x_\alpha\otimes t^n\longmapsto \Gamma_\alpha(n),
\qquad
d\longmapsto -D_0
\]
defines a representation of $\hat{\mathfrak g}$ on $V$ and yields an
isomorphism of representations of $\hat{\mathfrak{g}}$, 
\[
 V\simeq  \bigoplus_{\gamma } V_{\Lambda_\gamma} ,
\]
where $\gamma$ runs over minuscule weights (including $0$), and  
$$V_{\Lambda_\gamma}\simeq  S(t^{-1}\mathfrak{h}[t^{-1}])\otimes \mathbb{C}[\gamma+Q] $$
 is the level one highest weight integrable representation of $\hat{\mathfrak{g}}$ corresponding to the coset $\gamma+Q$.
\end{theorem}

For any $\lambda\in P$, set 
\begin{equation}\label{eq:deg_form}
 v_{\lambda}:= 1\otimes e^{\lambda} . \end{equation}
Then, the weight of $v_\lambda$ is given by
\begin{equation}\label{eq:wt-formula}
 {\rm wt}(v_\lambda)= \Lambda_0+\lambda- \frac{ (\lambda\mid \lambda)}{2} \delta .  \end{equation}
For $\gamma$ a minuscule weight in $P$, $v_{\gamma}$ is the highest weight vector of $V_{\Lambda_\gamma}$. 
In particular, $v_0$ is the highest weight vector of $V_{\Lambda_0}$. Note that by Formula \eqref{eq:wt-formula},  the degree of $v_{\gamma}$ is $-\frac{1}{2} (\gamma\mid \gamma)$. For any $\lambda\in P$, $v_\lambda$ is an extremal weight vector in $V_{\Lambda_\gamma}$, whenever $\lambda\in \gamma+Q$. 

From \eqref{eq:Gamma-alpha}, we obtain
\begin{equation} \label{eq:VOA-extrem-vec}
x_\alpha(z) v_\lambda=   \epsilon(\alpha,\lambda) z^{(\lambda\mid \alpha)}    \exp\!\left(\sum_{j\ge1}\frac{z^j}{j}\alpha(-j)\right) v_{\lambda+\alpha} . \end{equation}
In particular, 
\begin{equation} \label{eq:VOA-vec-0}
x_\alpha(z) v_0=   \exp\!\left(\sum_{j\ge1}\frac{z^j}{j}\alpha(-j)\right) v_{\alpha} . \end{equation}

\subsection{Sugawara operators and Demazure modules}
\label{subsec:sugawara-FK}
In this section we recall the Segal-Sugawara operators and examine their action on level one affine Demazure modules in the simply-laced setting.

The Segal-Sugawara field is given by
\begin{equation}
T(z)=\sum_a :x_a(z)x^a(z): =\sum_{n\in\mathbb Z}T_nz^{-n-2},
\end{equation}
where $x_a,x^a$ are dual bases of $\mathfrak g$ compatible with root decomposition.  The coefficients $T_n$ may be understood as operators in the restricted completion $U_{c}(\mathfrak{g}')$ of $U( \mathfrak{g}')$ where $\mathfrak{g}'=\mathfrak{g}[t,t^{-1}] \oplus \mathbb{C}K$. 
Going forward, we will be interested only in representations where $K$ acts as $1$, so we consider the quotient $U_c(\mathfrak{g}')/(K-1)$ which acts on these representations.  The level 1 Segal-Sugawara operators are
\begin{equation}\label{eq:sug}
L_n:=\frac{1}{2(1+h^{\vee})}T_n=\frac{1}{2(1+h^{\vee})} \sum_a \sum_{m \in \mathbb{Z}} :x_a(m)x^a(n-m): \end{equation}
where $h^\vee$ is the dual Coxeter number, and $: \,:$ denotes the normal ordering, defined by
\[
:x_a(m)x^a(n-m):
=
\begin{cases}
x_a(m)x^a(n-m), & m<0,\\[2mm]
x^a(n-m)x_a(m), & m\geq0.
\end{cases}
\]
 The operators $\{L_n\}_{n\in \mathbb{Z}}$ form a Virasoro algebra representation on $V_{\Lambda_\gamma}$, cf.\,\cite[Corollary 12.8]{Kac}.

By \cite[Lemma 12.8]{Kac},  
for any $n, j\in \mathbb{Z}$, and $x\in \g$,  in $U_c(\mathfrak{g}')/(K-1)$ we have formula:
\begin{equation}\label{eq:sug-com}
 [L_n, x(j)]= -jx(n+j).\end{equation}

For any $\lambda\in P^+$, we define the Demazure module $D(1,\lambda)$ as the following $\mathfrak{g}[t]$-module:
\begin{equation}  D(1,\lambda):= U(\mathfrak{g}[t])\cdot v_{-\lambda}.   \end{equation}
Then, $D(1,\lambda)$ is contained in $V_{\Lambda_\gamma}$, whenever $\lambda\in -\gamma+ Q$.

\begin{lemma}
\label{lem:truncation}
For a  dominant weight $\lambda \in P^+$, set
$a_\lambda=(\lambda\mid \theta)$. Then
\begin{equation}
\bigl(  t^{a_\lambda}\g [t]  \bigr)\cdot D(1,\lambda)=0.
\label{eq:truncation}
\end{equation}
\end{lemma}
\nopagebreak[4]
\begin{proof}
For each root $\alpha$, by \eqref{eq:VOA-extrem-vec},
\begin{equation}
x_\alpha(z)v_{-\lambda}=\epsilon(\alpha,-\lambda)\,z^{-(\alpha\mid \lambda) }
\exp\left(\sum_{j\geq1}\frac{\alpha(-j)}jz^j\right)(1\otimes e^{\alpha-\lambda}),
\label{eq:fk-pole}
\end{equation}
where  $x_\alpha(z)=\sum_{n\in \mathbb{Z} }x_\alpha(n)z^{-n-1}$.
All powers in the exponential are nonnegative. The coefficient of
$z^{-n-1}$ is therefore zero whenever
\begin{equation}
n \geq(\alpha\mid \lambda).
\label{eq:root-annihilation}
\end{equation}
Dominance gives $(\alpha\mid \lambda)  \leq(\theta\mid \lambda)=a_\lambda$ for every root
$\alpha$. Also, $h(n)v_{-\lambda}=0$ for every $h\in\fh$ and $n>0$.
Thus $I:=t^{a_\lambda}\g [t]$ annihilates $v_{-\lambda}$.

The subspace $I$ is a Lie ideal of $\g[t]$. Since $D(1,\lambda)$ is a cyclic module generated by $v_{-\lambda}$, we see $I$ annihilates $D(1,\lambda)$, proving \eqref{eq:truncation}.

\end{proof}

\begin{lemma}
\label{lem:primary}
For every $\lambda\in P$ and every $k\geq 1$,
\begin{equation}
L_k\cdot  v_\lambda=0.
\label{eq:positive-lattice}
\end{equation}
\end{lemma}
\begin{proof}
Let $\gamma$ (possibly $0$) be the minuscule weight such that $\lambda\in\gamma+Q$,
and let $w=t^{\lambda-\gamma}$ be the translation element in the affine Weyl group  associated to
$\lambda-\gamma\in Q$. By Formula \eqref{eq:wt-formula}, we have
$${\rm wt}(v_\lambda)=w\,{\rm wt}(v_\gamma)=w\Lambda_{\gamma},$$
 where $\Lambda_\gamma$ is the highest weight.
By \eqref{eq:sug}, $L_k$ has affine weight $k\delta$. Thus, if
$L_kv_\lambda\neq0$, then $w\Lambda_\gamma+k\delta$ is a weight of
$V_{\Lambda_\gamma}$. The weight set of an integrable representation is
invariant under the affine Weyl group, and $w^{-1}\delta=\delta$, see \cite[\S 6.5]{Kac}.
Hence $\Lambda_\gamma+k\delta$ would also be a weight, which is a contradiction when $k>0$. Therefore $L_kv_\lambda=0$.
\end{proof}

The following theorem describes a range of Segal-Sugawara operators which annihilate the affine Demazure modules $D(1, \lambda) \subset V$. Theorem \ref{thm:vir-bound} and Proposition \ref{prop:nonincl} which follows, are the main results of this section.
 
\begin{theorem}
\label{thm:vir-bound}
For any $\lambda\in P^+$ and  every integer $k\geq \max \{ (\lambda\mid \theta) -1, 1\}$,  we have 
\begin{equation}
L_k \cdot D(1,\lambda)=0.
\label{eq:vir-annihilation}
\end{equation}
\end{theorem}
\begin{proof}

Recall that we view our operators as elements in $U_c(\mathfrak{g}')/(K-1)$, where the following bracket computations hold. By \eqref{eq:sug-com}, when $j=0$,  $[L_k,x(0)]=0$; when $j\geq1$, we have
\[[L_k,x(j)]=-j\,x(k+j). \]

Hence, for any $j\geq 0$, $[L_k,x(j)]$ annihilates $D(1,\lambda)$ if $k+j\geq (\lambda\mid \theta) .$ by Lemma \ref{lem:truncation}. Then,  for any $ j_1,\ldots,j_N\geq0$, and $x_1, \cdots, x_{N}\in \g$,  we have
\begin{align*}
L_k \big (x_1(j_1)\cdots x_N(j_N)v_{-\lambda}   \big ) &=\sum_{s=1 }^N  x_1(j_1)  \cdots [L_k, x_s(j_s)]  \cdots x_N(j_N) v_{-\lambda} + x_1(j_1)\cdots x_N(j_N) L_k v_{-\lambda}   \\
&=x_1(j_1)\cdots x_N(j_N) L_k v_{-\lambda}  \\
&=0,
\end{align*}
where the last equality follows from Lemma \ref{lem:primary}. 
This proves the theorem, as $D(1, \lambda)= U(\g[t])v_{-\lambda}$. 
\end{proof}

Using Theorem \ref{thm:vir-bound}, we obtain the following useful criterion.
\begin{proposition} \label{prop:nonincl} Let $(\lambda,\mu)$ be a pair of dominant weights with $\mu\prec \lambda$. Assume that  $ (\lambda-\mu\mid \theta) \geq 1$. For any integer $q$, if $q\geq  (\lambda-\mu\mid \theta)$,  then
\[  x_\theta(-q)v_{-\mu} \notin  D(1,  \lambda)  .\]
\end{proposition}
\begin{proof}
For any integer $k\geq 1$, we have
\[ L_k(x_\theta(-q)v_{-\mu} )=[L_k, x_\theta(-q)]v_{-\mu}+ x_\theta(-q )L_kv_{-\mu} =q  x_\theta(k-q)v_{-\mu} . \]
Note that, $x_\theta(k-q)$ is a root vector associated to the affine root $\theta+ (k-q)\delta$, and $v_{-\mu}$ is an extremal weight vector of weight $\Lambda_0- \mu-\frac{ (\mu\mid \mu)}{2}\delta$.  
For the affine root $\theta+n\delta$, the associated affine coroot is
$(\theta+n\delta)^\vee=\theta^\vee+nK$, where $\theta^\vee$ is
identified with the Cartan element $\theta$ as above.
When $k \leq q-1+ (\mu\mid \theta)$, Formula \eqref{eq:wt-formula} gives
\[
{\rm wt}(v_{-\mu})\bigl(\theta^\vee+(k-q)K\bigr)
=k-q-(\mu\mid \theta)<0.
\]
By $sl_2$-string property on extremal weight vectors,  
$$L_k\cdot x_\theta(-q)v_{-\mu} = qx_\theta(k-q)v_{-\mu}\not=0.$$  

By assumption, $ (\lambda-\mu\mid \theta) \geq 1$. Since $\mu \preceq \lambda$ and both are dominant, one can easily see that we always have $(\lambda\mid \theta)-1 \geq 1$. 
When $k \geq (\lambda\mid \theta)-1$,  by Theorem \ref{thm:vir-bound}, 
  $$L_k\cdot D(1,\lambda)  =0.$$
 
    Thus, when there exists $k$ such that
    $$ (\lambda\mid \theta)-1\leq k \leq  q-1+ (\mu\mid \theta) ,$$
    we see $x_{\theta}(-q)v_{-\mu} \notin D(1,\lambda)$.  Equivalently, when $q\geq (\lambda-\mu\mid \theta)$,   we obtain  $x_\theta(-q)v_{-\mu} \not\in D(1,\lambda)$.
\end{proof}

\section{The $\theta$-minimality criterion for Finkelberg-Mirkovi\'c Schubert schemes}\label{sect:main_thm}

In this section, we first revisit Finkelberg-Mirkovi\'c Schubert schemes
for a general simple algebraic group via $z$-extensions. We then compare
the tangent spaces of spherical Schubert varieties and
Finkelberg-Mirkovi\'c Schubert schemes along highest-root directions at
$T$-fixed points. In Section~\ref{subsec:main-adjoint}, we prove
the $\theta$-minimality criterion for adjoint groups of every type,
using Section~\ref{sec:FK} and Polo's criterion in the simply-laced
case and diagram folding in the remaining cases.

\label{sec:fmtest}
\subsection{FM Schubert schemes  and $z$-extensions}
\label{subsec:fm-zextension}

Let $G$ be a  simple algebraic group over $\C$ with
Lie algebra $\g$, and let
$\pi:G_{\mathrm{sc}}\to G$ be its simply-connected cover.
Choose a maximal torus $T\subset G$ with Lie algebra
$\fh$ and a Borel subgroup $B\supset T$.
Let $\Phi$ be the root system of $(G,T)$, let $\Phi^+$ be the set of
positive roots determined by $B$, and let $I=\{1,\ldots,\ell\}$ index
the simple roots $\alpha_i$. We denote the simple coroots by
$\alpha_i^\vee$, and the fundamental weights and coweights by
$\omega_i$ and $\omega_i^\vee$, respectively. Let $\theta$ be the
highest root. Write $Q\subset P$ for the root and weight lattices, and
$Q^\vee\subset P^\vee$ for the coroot and coweight lattices; thus
\[
Q\subseteq X^*(T)\subseteq P,\qquad
Q^\vee\subseteq X_*(T)\subseteq P^\vee.
\]
Throughout this section and Section~\ref{sect:adm}, roots and coroots,
and weights and coweights, remain distinct.
Here $\langle\nu,\eta\rangle$ denotes the natural pairing of a
cocharacter $\nu$ with a character $\eta$, extended $\mathbb Q$-bilinearly
when needed. In particular, $\langle\alpha_i^\vee,\omega_j\rangle=\delta_{ij}$.
For $\lambda,\mu\in X_*(T)$, we write $\mu\preceq\lambda$ when
$\lambda-\mu$ is a nonnegative integral sum of simple coroots, and
$\mu\prec\lambda$ when in addition $\mu\neq\lambda$. For any $\eta\in X^*(T)^+$, we denote by $V_\eta$ the irreducible hightest-weight representation of $G$ with highest weight $\eta$.	

Let $LG$ denote the loop group functor and $L^+G$ denote the jet group functor.
The affine Grassmannian $\Gr_G$ is defined to be the fpqc sheaf quotient
$LG/L^+G$. For each $\lambda\in X_*(T)$, we have an associated point
$t^\lambda\in LG$ and $e_\lambda\in\Gr_G$.
Let $\Gr_{\lambda}=L^+G\cdot e_\lambda$, and let
$\Gr_{\preceq\lambda}$ denote its reduced closure, which is called a
spherical Schubert variety.

For any $\C$-algebra $R$, set
$D_R:={\rm Spec}R[[t]]$ and $D_R^*:={\rm Spec}R((t))$.
Then the affine Grassmannian $\Gr_G$ can be described as the following moduli functor:
\[
\Gr_G(R)=\{(\mathscr E,\phi)\mid \mathscr E\text{ is a }G\text{-torsor over }D_R,
\ \phi:\mathscr E|_{D_R^*}\simeq\mathring{\mathscr E}_{D_R^*}\},
\]
where $\mathring{\mathscr E}_{D_R^*}$ denotes the trivial $G$-torsor over $D_R^*$.

To formulate the Finkelberg-Mirkovi\'c moduli problem uniformly for an
arbitrary isogeny type, and in particular to have access to all
fundamental representations of the simply-connected cover, we use a
$z$-extension of $G$, following \cite{ZhuPCMI,KPZ}.

\begin{definition}
A \emph{$z$-extension} of $G$  is a central extension of $G$ 
\begin{equation}\label{eq:zfm-extension}
1\to S\to G'\xrightarrow{p}G\to 1
\end{equation}
by a torus $S$,  such that the derived group $(G')^{\mathrm{der}}$ is the simply-connected cover  $G_{\mathrm{sc}}$ of $G$.
\end{definition}
Such extensions always exist; see \cite[Proposition~3.1]{MS} and \cite[\S4.2.4]{KPZ}. For example, when $G=PGL_n$,  \[1\to Z(GL_n)\to GL_n\to PGL_n\to 1\] is a $z$-extension of $PGL_n$. If $G$ is simply-connected, then we can take $G'=G$ and $S=1$. 

Set $T'=p^{-1}(T)$. The exact sequence of  tori gives an exact sequence of lattices: 
\[
0\to X_*(S)\to X_*(T')
\xrightarrow{p_*}X_*(T)\to 0.
\]
Hence any $\lambda\in X_*(T)^+$ has a dominant lift
$\lambda'\in X_*(T')^+$. Let $A$ be the quotient $G'/G_{\mathrm{sc}}$, which is a torus.  We have a surjective map of tori $T'\to A$ and the induced morphism $X_*(T')\to X_*(A)$. Let $e_\lambda$ be point in $\mathrm{Gr}_G$ associated to $\lambda$.  We denote $\mathrm{Gr}_G^{[\lambda]}$ the connected component of $\mathrm{Gr}_G$ containing $e_{\lambda}$.  

Fix a lift $\lambda'\in X_*(T')$ of $\lambda$. Let $\nu\in X_*(A)$ denote the image of $\lambda'$ via the map $X_*(T')\to X_*(A)$. Let $F_\nu$ be the schematic fiber of the morphism $\Gr_{G'}\to \Gr_A$ over $e_\nu$.  \begin{lemma}\label{lem:comp}
We have the following isomorphisms of ind-schemes:
\[
\Gr_{G_{\mathrm{sc}}}
\xrightarrow[\sim]{\,y\mapsto t^{\lambda'}y\,}   F_\nu  \xrightarrow{\sim}
\Gr_G^{[\lambda]},
\]
where the second map is induced from $p:G'\to G$.
\end{lemma}

\begin{proof}
Consider the following exact sequence
\[
1\to G_{\mathrm{sc}}\to G'\xrightarrow{q}A\to1 .
\]
Let $F_0$ be the fiber of $\Gr_{G'}\to\Gr_A$ over the base point
$e_0$. Then  $F_0$  can be identified with $\Gr_{G_{\mathrm{sc}}}$.
 Indeed, for an $R$-point $(\mathscr E,\phi)$ of this fiber, the induced $A$-torsor $\mathscr E\times^{G'}A$, together with its trivialization over the punctured disc $D^*_R$, is the base point of $\Gr_A(R)$. Hence its
trivialization over $D_R$ gives a reduction of $\mathscr E$ to $\ker(q)=G_{\mathrm{sc}}$, and conversely extension of structure group
from $G_{\mathrm{sc}}$ to $G'$ produces an $R$-point of the fiber $F_0$. Thus $F_0\simeq\Gr_{G_{\mathrm{sc}}}$.
Since the image of $t^{\lambda'}$ in the loop group $LA$ is $t^\nu$, translation by
$t^{\lambda'}$ carries $F_0$ isomorphically onto $F_\nu$. Thus,  $\Gr_{G_{\mathrm{sc}}}\simeq F_\nu$.

For the simply-connected covering
$G_{\mathrm{sc}}\to G$, we have an isomorphism $\pi_{\Gr}:\Gr_{G_{\mathrm{sc}}}  \xrightarrow{\sim} \Gr_G^0$, 
where  $\Gr_G^0$ is the neutral component of $\Gr_G$, cf.\,\cite[Theorem 1.3.11(3)]{ZhuPCMI}. Then, we have isomorphisms
\begin{equation}\label{eq:j_lambda}
\Gr_{G_{\mathrm{sc}}} \xrightarrow{\pi_{\Gr}} \Gr_G^0
\xrightarrow{\ t^\lambda\ }\Gr_G^{[\lambda]}.
\end{equation}
Since $p(t^{\lambda'})=t^\lambda$,  the morphism induced by $p$ gives rise to the isomorphism $F_\nu\simeq \Gr_{G}^{[\lambda]}$.
\end{proof}


Finkelberg and Mirkovi\'c \cite{FM} proposed a moduli interpretation of
spherical Schubert varieties by imposing pole conditions on trivializations of vector bundles associated to irreducible
representations. We apply this construction to the $z$-extension $G'$ of $G$. For a $G'$-torsor $\mathscr E$ over $D_R$ and a dominant
weight $\eta'\in X^*(T')^+$, write
$\mathscr E(\eta')=\mathscr E\times^{G'}V_{\eta'}$. For a point
$(\mathscr E,\phi)\in\Gr_{G'}(R)$, let $L_{\eta',R}$ denote the lattice
\[
\phi\bigl(\Gamma(D_R,\mathscr E(\eta'))\bigr)
\subset V_{\eta'}\otimes R((t)),
\]
and set $L^0_{\eta',R}=V_{\eta'}\otimes R[[t]]$.
Following \cite[Section~4.1]{BHY}, define
\begin{equation}\label{eq:zfm-upstairs}
\Gr^{\mathrm{FM},G'}_{\preceq\lambda'}(R)
:=
\left\{(\mathscr E,\phi)\in\Gr_{G'}(R)\;\middle|\;
L^0_{\eta',R}\subseteq
 t^{-\langle\lambda',\eta'\rangle}L_{\eta',R}
\text{ for all }\eta'\in X^*(T')^+\right\}.
\end{equation}
Our convention differs from some references by a $-w_0$ twist, e.g.\,\cite[\S2.3]{KMW}, \cite[Remark~2.1.7]{ZhuPCMI}; see the discussion in \cite[Remark~4.2]{BHY}. 

Applying the defining condition to the dual representation
$V_{\eta'}^\vee\simeq V_{-w_0(\eta')}$ gives the opposite inclusion. Thus
\eqref{eq:zfm-upstairs} is equivalent to imposing the two-sided condition
\begin{equation}\label{eq:zfm-two-sided}
t^{-\langle\lambda',w_0(\eta')\rangle}L_{\eta',R}
\subseteq L^0_{\eta',R}
\subseteq t^{-\langle\lambda',\eta'\rangle}L_{\eta',R}
\end{equation}
for all $\eta' \in X^*(T')^+$.
For $\lambda' \in X_*(T')^+$ a dominant lift of $\lambda$, we will use $\mathrm{Gr}_{\preceq \lambda'}^{\mathrm{FM}, G'}$ in order to define FM Schubert schemes associated to $\lambda$. First we check that $\Gr_{\preceq \lambda'}^{\mathrm{FM}, G'}$ lies in the correct component.

\begin{lemma}\label{lem:zfm-fiber}
The scheme $\Gr^{\mathrm{FM},G'}_{\preceq\lambda'}$ is a closed subscheme of $F_\nu$. 
\end{lemma}
\begin{proof}
Let $q:G'\to A$ be the quotient map, and let
$x\in\Gr_{\preceq\lambda'}^{\mathrm{FM},G'}(R)$.
For $\chi\in X^*(A)$, regard $\chi$ as a character of $G'$ via $q$.
Since $\chi$ is trivial on $G_{\mathrm{sc}}$, it is fixed by the Weyl
group of $G'$; in particular, $w_0(\chi)=\chi$.
Moreover, both $\chi$ and $-\chi$ are dominant.

Let $L_{\chi,R}\subset R((t))$ be the lattice associated with the
one-dimensional representation of $G'$ corresponding to the character $\chi$. Applying the
two-sided FM condition \eqref{eq:zfm-two-sided} to $\chi$ gives
\[
t^{-\langle\lambda',\chi\rangle}L_{\chi,R}
\subseteq R[[t]]
\subseteq
t^{-\langle\lambda',\chi\rangle}L_{\chi,R}.
\]
Hence, 
$L_{\chi,R} = t^{\langle\lambda',\chi\rangle}R[[t]] = t^{\langle\nu,\chi\rangle}R[[t]]$.  Since this holds for every $\chi\in X^*(A)$, the image of $x$ in
$\Gr_A(R)$ is $e_{\nu,R}$. Thus the morphism  $\Gr_{\preceq\lambda'}^{\mathrm{FM},G'}\to \Gr_A$
factors through $e_\nu$, which proves  $\Gr_{\preceq\lambda'}^{\mathrm{FM},G'}\subseteq F_\nu$ scheme-theoretically.
\end{proof}
%

\begin{lemma}\label{lem:FM_red}
The morphism induced by $p:G'\to G$ restricts to an isomorphism
\[
\bigl(\Gr^{\mathrm{FM},G'}_{\preceq\lambda'}\bigr)_{\mathrm{red}} \simeq  \Gr_{\preceq\lambda},
\]
where $\bigl(\Gr^{\mathrm{FM},G'}_{\preceq\lambda'}\bigr)_{\mathrm{red}}$ denotes the underlying reduced subscheme of
$\Gr^{\mathrm{FM},G'}_{\preceq\lambda'}$.
\end{lemma}

\begin{proof}
By \cite[Propositions~2.1.4--2.1.5 and Remark~2.1.7]{ZhuPCMI}, $\bigl(\Gr^{\mathrm{FM},G'}_{\preceq\lambda'}\bigr)_{\mathrm{red}} =\Gr^{G'}_{\preceq\lambda'}$.
The assertion then follows from Lemma~\ref{lem:comp}.
\end{proof}

Next we show that the image of $\Gr_{\preceq \lambda'}^{\mathrm{FM}, G'}$ in $\Gr_G$ induced by $p: G' \rightarrow G$ is independent of the $z$-extension $G'$ and the lift $\lambda'$. We will make auxiliary use of the following lattices.
For $\eta\in P^+$, let $V_\eta$ denote the irreducible
$G_{\mathrm{sc}}$-module of highest weight $\eta$. For a $\C$-algebra $R$, define
\begin{equation}\label{eq:zfm-intrinsic-lattice}
\mathcal B_{\lambda,\eta,R}
:=
\bigoplus_{\beta\in{\rm wt}(V_\eta)}
 t^{\langle\lambda,\eta-\beta\rangle}
 V_\eta(\beta)\otimes R[[t]].
\end{equation}
The exponents are integers because $\eta-\beta$ lies in the root lattice
$Q\subset X^*(T)$. Thus this lattice is defined for every
$\lambda\in X_*(T)$, without choosing a lift of $\lambda$ or $\eta$.

\begin{proposition}\label{prop:zfm-independence}
The morphism $\Gr^{\mathrm{FM},G'}_{\preceq\lambda'}\to\Gr_G$
induced by $p:G'\to G$ is a closed immersion. Its image is independent of
both the choice of $z$-extension \eqref{eq:zfm-extension} and the lift
$\lambda'$.
\end{proposition}

\begin{proof}
In view of the pole condition in \eqref{eq:zfm-upstairs} and Lemma~\ref{lem:zfm-fiber}, 
$\Gr^{\mathrm{FM},G'}_{\preceq\lambda'}$ is a closed subscheme of $F_\nu$. 
By Lemma \ref{lem:comp}, the morphism $p$ induces an identification between $F_\nu$ and the component
$\Gr_G^{[\lambda]}$. This proves the first assertion.

From  Lemma~\ref{lem:comp}, 
let $f_{\lambda'}:F_\nu\simeq \Gr_{G_{\mathrm{sc}}}$
be the inverse of the map $\Gr_{G_{\mathrm{sc}}}\to F_\nu$ given by $y\mapsto t^{\lambda'}y$. 
The morphism $f_{\lambda'}$ gives the following moduli interpretation.  For any $x=(\mathscr E,\phi)\in F_\nu(R)$, we first replace
$\phi$ by $t^{-\lambda'}\phi$.  The induced $A$-trivialization extends
to $D_R$. Thus,  $(\mathscr E,t^{-\lambda'}\phi)$ determines a $G_{\mathrm{sc}}$-reduction, which gives rise to a $R$-point in $\Gr_{G_{\mathrm{sc}}}$.

Let $j_{\lambda}=t^{\lambda} \circ \pi_{\Gr}$ be the composition from \eqref{eq:j_lambda}. Since $p(t^{\lambda'})=t^\lambda$, the map $F_\nu\to\Gr_G$ induced by $p$ is $j_\lambda\circ f_{\lambda'}$.  As $j_\lambda$ depends only on
$\lambda$, it suffices to show that the closed subscheme
$f_{\lambda'}(\Gr^{\mathrm{FM},G'}_{\preceq\lambda'})$ of
$\Gr_{G_{\mathrm{sc}}}$ is independent of the $z$-extension and the lift.

We identify $(G')^{\mathrm{der}}$ with $G_{\mathrm{sc}}$, and write
$T_{\mathrm{sc}}=T'\cap G_{\mathrm{sc}}$. For any dominant weight
$\eta\in P^+=X^*(T_{\mathrm{sc}})^+$, the exact sequence
\[
0\to X^*(A)\to X^*(T')
\to X^*(T_{\mathrm{sc}})\to 0
\]
provides a lift $\eta'\in X^*(T')$ of $\eta$. This lift remains dominant, and
$V_{\eta'}|_{G_{\mathrm{sc}}}\simeq V_\eta$. Conversely, every irreducible
representation of $G'$ restricts irreducibly to $G_{\mathrm{sc}}$.

For $x\in F_\nu(R)$, let $y=f_{\lambda'}(x)$, and let
$L_{\eta,R}(y)\subset V_\eta\otimes R((t))$ be its associated lattice.
As usual we write $L^0_{\eta,R}=V_\eta\otimes R[[t]]$. After identifying
$V_{\eta'}|_{G_{\mathrm{sc}}}\simeq V_\eta$, the definition of
$f_{\lambda'}$ gives
\begin{equation}\label{eq:zfm-translated-lattice}
L_{\eta,R}(y)
=
t^{-\lambda'}L_{\eta',R}(x).
\end{equation}
As a consequence, $x$ satisfies the FM condition \eqref{eq:zfm-upstairs}
if and only if
\begin{equation}\label{eq:new-latt}
t^{\langle\lambda',\eta'\rangle}
(t^{-\lambda'}L^0_{\eta,R})
\subseteq L_{\eta,R}(y),
\qquad(\eta'\in X^*(T')^+).
\end{equation}
If $\beta'$ is a  $T'$-weight appearing in $V_{\eta}(\beta)$, then
$\eta'-\beta'=p^*(\eta-\beta)$ is in the root lattice $Q\subset X^*(T')$. Hence
\[
\langle\lambda',\eta'-\beta'\rangle
=\langle\lambda,\eta-\beta\rangle\in \mathbb{Z},
\]
so the lattice on the left of \eqref{eq:new-latt} is exactly $\mathcal B_{\lambda,\eta,R}$,which is
independent of the choice of the lift $\eta'$. Since every dominant $\eta$ lifts to $\eta'$
and $f_{\lambda'}$ is an isomorphism, a point
$y\in\Gr_{G_{\mathrm{sc}}}(R)$ belongs to
$f_{\lambda'}(\Gr^{\mathrm{FM},G'}_{\preceq\lambda'})(R)$ if and only if
\begin{equation}\label{eq:zfm-intrinsic-test}
\mathcal B_{\lambda,\eta,R}
\subseteq L_{\eta,R}(y),
\qquad(\eta\in P^+).
\end{equation}

These conditions depend only on $\lambda$ and $y$, and define the same
subfunctor of $\Gr_{G_{\mathrm{sc}}}$ for every $z$-extension and every
lift $\lambda'$. Since the map $j_\lambda$ is also independent of these
choices, the images in $\Gr_G$ agree on every $\C$-algebra, and hence
agree as closed subschemes. This proves the second assertion.
\end{proof}

From the proof of Proposition~\ref{prop:zfm-independence}, we have the following intrinsic description of $j_\lambda^{-1}( \Gr^{\mathrm{FM}}_{\preceq\lambda})$ as a closed subscheme in $\Gr_{G_{sc}}$.
\begin{corollary}\label{cor:fm-intrinsic}
For any $\C$-algebra $R$, a point $y\in\Gr_{G_{\mathrm{sc}}}(R)$
belongs to $j_\lambda^{-1}(\Gr^{\mathrm{FM}}_{\preceq\lambda})(R)$
if and only if
\[
\mathcal B_{\lambda,\eta,R}\subseteq L_{\eta,R}(y)
\qquad(\eta\in P^+),
\]
where $L_{\eta,R}(y)\subset V_\eta\otimes R((t))$ is the lattice associated to $y$.
\end{corollary}

\begin{definition}
For any simple $G$ and $\lambda\in X_*(T)^+$, we define
$\Gr^{\mathrm{FM}}_{\preceq\lambda}\subset\Gr_G$ to be the schematic
image of
$\Gr^{\mathrm{FM},G'}_{\preceq\lambda'}\to\Gr_G$.  
\end{definition}
By Proposition~\ref{prop:zfm-independence} and Lemma~\ref{lem:FM_red}, we have
\begin{equation}\label{eq:zfm-scheme-comparison}
\Gr^{\mathrm{FM},G'}_{\preceq\lambda'}
\simeq  \Gr^{\mathrm{FM}}_{\preceq\lambda},
\qquad
\bigl(\Gr^{\mathrm{FM}}_{\preceq\lambda}\bigr)_{\mathrm{red}}
\simeq \Gr_{\preceq\lambda}.
\end{equation}
If $\lambda\in Q^\vee\subset X_*(T)$, then this definition agrees with the original Finkelberg-Mirkovi\'c definition of $\Gr_{\preceq \lambda}^{\rm{FM}}$ associated to $\lambda$ for $G_{\mathrm{sc}}$  via the isomorphism
$\Gr_{G_{\mathrm{sc}}}\simeq\Gr_G^0$.

We have the following simple corollary:
\begin{corollary}\label{cor:zfm-tangent}
Let $\mu \preceq \lambda$ in $X_*(T)^+$, and let $\lambda',\mu'$ be
compatible lifts to $G'$, i.e.
\begin{equation}\label{eq:zfm-compatible-lift}
\lambda'-\mu'=\lambda-\mu\in Q^\vee,
\end{equation}
where $Q^\vee=X_*(T_{\mathrm{sc}})$ is viewed inside $X_*(T')$.
 Then the natural maps induce isomorphisms
\[
T_{e_{\mu'}}
\Gr_{\preceq \lambda'}^{\mathrm{FM},G'}
\simeq
T_{e_\mu}
\Gr_{\preceq \lambda}^{\mathrm{FM}},
\qquad
T_{e_{\mu'}}
\Gr_{\preceq \lambda'}^{G'}
\simeq
T_{e_\mu}
\Gr_{\preceq \lambda}.
\]
\end{corollary}

\begin{proof}
The first is immediate from Proposition~\ref{prop:zfm-independence}, and
the second is the standard $z$-extension reduction; see
\cite[\S4.2.4]{KPZ}.
\end{proof}


\subsection{Tangent spaces of FM Schubert schemes} 
\label{subsec:fm-tangents}
Let $G$ be a simple algebraic group as in the previous subsection. In this subsection we work with the tangent spaces $T_{e_{\mu}}\rm{Gr}_{\preceq \lambda}^{\rm{FM}}$, so in this subsection we always take $R= \mathbb{C}[\epsilon]/\epsilon^2$ to be the ring of dual numbers.

For any $\mu\in X_*(T)$, let $T_{e_\mu}\Gr_G$ denote the tangent space of
$\Gr_G$ at $e_\mu$. Then we have a natural isomorphism
\begin{equation}\label{eq:tangent_mu}
T_{e_\mu}\Gr_G
\simeq {\rm Ad}_{t^\mu}(t^{-1}\g[t^{-1}])
\simeq t^{-1}\mathfrak h[t^{-1}]
\oplus\bigoplus_{\alpha\in\Phi,\,k>0}
\g_\alpha\otimes t^{\la\mu,\alpha\ra-k}.
\end{equation}


For $q \geq 1$, consider the $R$-point of $LG$ given by \[1+x_{\theta}(-q) \epsilon\]
and, via the action of $LG$ on $\Gr_G$, the $R$-point of $\mathrm{Gr}_G$ given by 
\begin{equation}\label{eq:infinitesimal}
(1+x_{\theta}(-q) \epsilon) e_{\mu}.
\end{equation}
We write $\xi_{\mu, q}$ for the associated tangent vector in $T_{e_{\mu}} \Gr_G$.
By \eqref{eq:tangent_mu}, the vector
$\xi_{\mu,q}\neq0$ if and only if
\begin{equation}\label{eq:nonzero}
{\rm Ad}_{t^{-\mu}}(x_\theta(-q))
=x_\theta t^{-q-\la\mu,\theta\ra}\in t^{-1}\g[t^{-1}].
\end{equation}

The following lemma generalizes
\cite[Lemma~4.3]{BHY}.

\begin{lemma}\label{lem:Plucker}
Let $\mu\preceq\lambda$ in $X_*(T)^+$, and choose compatible lifts $\lambda',\mu'$ as in
\eqref{eq:zfm-compatible-lift}. Under the identification
\eqref{eq:tangent_mu}, the tangent space
$T_{e_\mu}\Gr^{\mathrm{FM}}_{\preceq\lambda}$ is the subspace of
${\rm Ad}_{t^\mu}(t^{-1}\g[t^{-1}])$ consisting of those $X$ such
that, for every $\eta'\in X^*(T')^+$,
\[
L^0_{\eta',R}\subseteq
 t^{-\langle\lambda',\eta'\rangle}
 (1+X\epsilon)t^{\mu'}L^0_{\eta',R},
\]
or equivalently,
\[
(1-X\epsilon)L^0_{\eta',R}\subseteq
 t^{-\langle\lambda',\eta'\rangle}t^{\mu'}L^0_{\eta',R}.
\]
Here $X$ acts through the derived Lie algebra $\g\subset{\rm Lie}(G')$.
\end{lemma}
\begin{proof}
 Via the isomorphism $F_\nu\simeq \Gr_G^{[\lambda]}$,  the identification \eqref{eq:tangent_mu} gives rise to 
\[ 
T_{e_{\mu'}} F_\nu \simeq {\rm Ad}_{t^{\mu}}
\bigl(t^{-1}\g[t^{-1}]\bigr),\]
where the tangent vector corresponding to
$X\in {\rm Ad}_{t^{\mu}} \bigl(t^{-1}\g[t^{-1}]\bigr)$ is represented by the $R$-point
$y_X:=[(1+\epsilon X)t^{\mu'}]\in \Gr_{G'}(R)$.  Modulo $\epsilon$, it gives the point $e_{\mu'}$.

For every dominant weight $\eta'\in X^*(T')^+$, the lattice associated
to $y_X$ in  $V_{\eta'}\otimes R((t))$ is
\[
L_{\eta',R}(y_X)
=
(1+\epsilon X)t^{\mu'}L^0_{\eta',R}.
\]
By the definition of the FM Schubert scheme, $y_X\in \Gr_{\preceq\lambda'}^{\mathrm{FM},G'}(R)$
if and only if, for every $\eta'\in X^*(T')^+$,
\begin{equation}\label{eq:incl_1}
L^0_{\eta',R} \subseteq  t^{-\langle\lambda',\eta'\rangle}
(1+\epsilon X)t^{\mu'}L^0_{\eta',R}.
\end{equation}
Thus, using Corollary~\ref{cor:zfm-tangent}, these inclusions are
equivalent to $X\in T_{e_\mu}\Gr_{\preceq\lambda}^{\mathrm{FM}}$.

Finally, since $\epsilon^2=0$, we have $(1+\epsilon X)^{-1}=1-\epsilon X$. Hence the inclusion \eqref{eq:incl_1} is equivalent to
\[
(1-\epsilon X)L^0_{\eta',R}
\subseteq
t^{-\langle\lambda',\eta'\rangle}
t^{\mu'}L^0_{\eta',R}.
\]
\end{proof}

Let $\mu,\lambda\in X_*(T)^+$ and assume $\mu\preceq\lambda$, so
$e_\mu\in\Gr_{\preceq\lambda}$. Write
$\lambda-\mu=\sum_i a_i\alpha_i^\vee$, with $a_i\in\Z^{\geq0}$, and set
\[
m_{\lambda-\mu}:=\min_i a_i
=\min_i\{\langle\lambda-\mu,\omega_i\rangle\},
\]
where $\omega_i$ are the fundamental weights of $G_{\mathrm{sc}}$.

\begin{proposition}\label{prop:fmtest}
For any $q\geq1$,
$\xi_{\mu,q}\in T_{e_\mu}\Gr^{\rm FM}_{\preceq\lambda}$ if and only if
$q\leq m_{\lambda-\mu}$.
\end{proposition}

\begin{proof}
First choose compatible lifts
$\lambda',\mu'$. Suppose
$\xi_{\mu,q}\in T_{e_\mu}\Gr^{\rm FM}_{\preceq\lambda}$.
By Lemma~\ref{lem:Plucker}, for every $\eta'\in X^*(T')^+$ we have
\begin{equation}\label{eq:lattice_incl}
(1-x_\theta(-q)\epsilon)L^0_{\eta',R}
\subseteq t^{-\langle\lambda',\eta'\rangle}t^{\mu'}L^0_{\eta',R}.
\end{equation}
If $V_{\eta'}$ is one-dimensional, then $x_\theta$ acts trivially and this
condition imposes no restriction. Otherwise let $V_\eta$ be the
$G_{\mathrm{sc}}$-module obtained by restriction, with $\eta\in P^+\setminus\{0\}$.
For a $T'$-weight vector $u_{\nu'}\in V_{\eta'}$ with
$x_\theta u_{\nu'}\neq0$, the coefficient in $\epsilon$ of
$(1-x_\theta(-q)\epsilon)u_{\nu'}$ is $-x_\theta(-q)u_{\nu'}$, which has
degree $-q$ and weight $\nu'+\theta$. Thus \eqref{eq:lattice_incl} is
equivalent to
\begin{equation}\label{eq:fm-weight}
q+\langle\mu',\nu'+\theta\rangle
\leq\langle\lambda',\eta'\rangle
\end{equation}
for every $\nu'$ such that $x_\theta V_{\eta'}(\nu')\neq0$.

We have
\begin{equation}\label{eq:max-weight}
\max_{\nu':\,x_\theta V_{\eta'}(\nu')\neq0}
\langle\mu',\nu'+\theta\rangle
=\langle\mu',\eta'\rangle.
\end{equation}
Indeed, all weights of $V_{\eta'}$ have the same central character, and
$\mu$ is dominant on the derived group, so the right-hand side bounds the
pairing of every weight. As in the proof of \cite[Theorem~4.5]{BHY}, for a highest weight vector $v_{\eta'}$, we have
$x_{-\theta}v_{\eta'}\neq 0$. Since
$[x_\theta,x_{-\theta}]=-\theta^\vee$, we have
\[
x_\theta x_{-\theta}v_{\eta'}
=-\langle\theta^\vee,\eta\rangle v_{\eta'}\neq0,
\]
which gives equality in \eqref{eq:max-weight}. Thus the condition \eqref{eq:fm-weight} is equivalent to
\begin{equation}\label{eq:fm-eta}
q\leq\langle\lambda'-\mu',\eta'\rangle
=\langle\lambda-\mu,\eta\rangle.
\end{equation}
If $q \leq \langle \lambda-\mu, \eta \rangle$ holds for all dominant $\eta$, we easily obtain $q \leq \langle \lambda-\mu, \omega_i \rangle$ for fundamental weights $\omega_i$ and thus 
\[q \leq m_{\lambda-\mu}. \] 

In the other direction, we assume $q \leq m_{\lambda-\mu}$. If all fundamental pairings are at least $q$, for an arbitrary
$\eta=\sum_i n_i\omega_i\neq0$ we have
\[q \leq \sum_i n_i \langle \lambda-\mu, \omega_i \rangle = \langle \lambda - \mu, \eta\rangle.\]
Hence \eqref{eq:lattice_incl} holds for every irreducible representation
of $G'$, including the one-dimensional representations, which shows $\xi_{\mu,q} \in T_{e_{\mu}} \Gr_{\preceq \lambda}^{\rm{FM}}$. 
\end{proof}

\subsection{Main theorem for adjoint groups}
\label{subsec:main-adjoint}

In this subsection, we assume that $G$ is adjoint, and
keep the notation of Section~\ref{subsec:fm-zextension}. Here
$\langle\, ,\,\rangle$ denotes the natural pairing, and
 $(\,,\mid\,,)$ denote the normalized invariant form. It induces an injection
\[
\nu\longmapsto\nu^\flat:X_*(T)\hookrightarrow P,
\qquad (\nu^\flat\mid\eta)=\langle\nu,\eta\rangle\quad(\eta\in P).
\]
In particular,
\[
(\alpha_i^\vee)^\flat=\frac{2\alpha_i}{(\alpha_i\mid\alpha_i)},
\qquad
(\omega_i^\vee)^\flat=\frac{2\omega_i}{(\alpha_i\mid\alpha_i)}.
\]
For $\lambda\in X_*(T)^+$, let $\gamma$ be the unique minuscule coweight (or $0$)
such that $\lambda+\gamma\in Q^\vee$, and let $\Gr_G^\gamma$ be the connected
component containing $e_{-\gamma}$ (and hence $e_\lambda$). Let
$\mathscr L$ be the level one line bundle on $\Gr_G$, whose restriction to each connected component is the positive generator of its Picard group. Write
$V_{\Lambda_{\gamma^\flat}}$ for the level-one integrable highest-weight
module with finite highest weight $\gamma^\flat$, using the degree
normalization of Section~\ref{sec:FK} when $G$ is simply-laced.
For $\nu+\gamma\in Q^\vee$, let $v_{-\nu^\flat}$ be an extremal vector whose line corresponds to $e_\nu$.
By the affine Borel-Weil theorem,
\begin{equation}\label{eq:pluck}
V_{\Lambda_{\gamma^\flat}}\simeq H^0(\Gr_G^\gamma,\mathscr L)^*,
\end{equation}
while
\[
D(1,\lambda)\simeq H^0(\Gr_{\preceq\lambda},\mathscr L)^*,
\qquad
D(1,\lambda):=U(\g[t])v_{-\lambda^\flat}\subset V_{\Lambda_{\gamma^\flat}};
\]
see \cite[Proposition~1.4.4]{Zhu09} and \cite{Ku}. These identifications
give the Pl\"ucker embeddings
\begin{equation}\label{eq:pluck_embed}
\iota:\Gr_G^\gamma\xhookrightarrow{}\mathbb P(V_{\Lambda_{\gamma^\flat}}),
\qquad
\Gr_{\preceq\lambda}\xhookrightarrow{}\mathbb P(D(1,\lambda)),
\end{equation}
with $e_\mu\mapsto[v_{-\mu^\flat}]$ for every $\mu\preceq\lambda$.

The following theorem is Polo's criterion for the tangent space of Schubert
varieties; see finite type in \cite{Po}, Kac-Moody type in
\cite[Exercise~12.1.E(9)]{Ku}. The formulation below follows from
\cite[Theorem~4.2]{HLR}.

\begin{theorem}\label{thm:polo}
For any $\mu,\lambda\in X_*(T)^+$ with $\mu\preceq\lambda$, we have an equality
\[
T_{e_\mu}\Gr_{\preceq\lambda}
=\{X\in T_{e_\mu}\Gr_G\mid X\cdot v_{-\mu^\flat}\in D(1,\lambda)\}.
\]
\end{theorem}

\begin{proof}
By \eqref{eq:pluck_embed}, we have the following commutative diagram:
\[
\begin{tikzcd}
\Gr_{\preceq\lambda}\arrow[r,hook,"\iota"]\arrow[d,hook]
&\mathbb P(D(1,\lambda))\arrow[d,hook]\\
\Gr_G^{\gamma}\arrow[r,hook,"\iota"]
&\mathbb P( V_{\Lambda_{\gamma^\flat}}).
\end{tikzcd}
\]
Note that $\Gr_G^{\gamma}$ is a partial affine flag variety of the loop group $LG_{\rm sc}$. By \cite[Theorem~4.2]{HLR},  $\Gr_{\preceq\lambda}$ is cut out inside
$\Gr_G^{\gamma}$ by linear equations via the Pl\"ucker embedding.
Hence, scheme-theoretically,
\begin{equation}\label{eq:sche_int}
\Gr_{\preceq\lambda}
=\Gr_G^{\gamma}\cap\mathbb P(D(1,\lambda)).
\end{equation}
Recall that
\[
T_{[v_{-\mu^\flat}]}\mathbb P(D(1,\lambda))
=
{\rm Hom}\bigl(\C v_{-\mu^\flat},
D(1,\lambda)/\C v_{-\mu^\flat}\bigr).
\]
By \eqref{eq:sche_int},
\[
T_{e_\mu}\Gr_{\preceq\lambda}
=T_{e_\mu}\Gr_G\cap
T_{[v_{-\mu^\flat}]}\mathbb P(D(1,\lambda)).
\]
By the same argument as in the proof of \cite[Corollary~4.3]{HLR},
\begin{equation}\label{eq:equiv_Dem}
X\in T_{e_\mu}\Gr_{\preceq\lambda}
\iff
Xv_{-\mu^\flat}\bmod\C v_{-\mu^\flat}
\in D(1,\lambda)/\C v_{-\mu^\flat},
\end{equation}
where we identify
$T_{e_\mu}\Gr_G\simeq{\rm Ad}_{t^\mu}(t^{-1}\g[t^{-1}])$.
The right-hand side of \eqref{eq:equiv_Dem} is equivalent to
$Xv_{-\mu^\flat}\in D(1,\lambda)$. This proves the theorem.
\end{proof}

When $G$ is simply-laced, the map $\nu\mapsto\nu^\flat$ is an
isomorphism with $(\alpha_i^\vee)^\flat=\alpha_i$ and
$(\omega_i^\vee)^\flat=\omega_i$. Under this identification, $D(1,\lambda)$
is the module denoted by $D(1,\lambda^\flat)$ in Section~\ref{sec:FK}.
Theorem~\ref{thm:FKS} then realizes the modules and extremal vectors above
in the Frenkel-Kac construction.

\begin{corollary}\label{cor:polo}
Assume that $G$ is simply-laced.
For  any $\lambda, \mu\in X_*(T)^+$ with $\mu \prec \lambda$, 
if $1\leq\la\lambda-\mu,\theta\ra\leq q$, then
$$\xi_{\mu,q}\notin T_{e_\mu}\Gr_{\preceq\lambda}.$$
\end{corollary}

\begin{proof}
The corollary follows from Proposition~\ref{prop:nonincl} and
Theorem~\ref{thm:polo}.
\end{proof}

\begin{definition}\label{def:admissible}
For any $\lambda,\mu\in X_*(T)^+$, we say the pair $(\lambda,\mu)$ is
$\theta$-minimal if $\mu\prec\lambda$ and
\begin{equation}\label{eq:question}
1\leq\langle\lambda-\mu,\theta\rangle
\leq m_{\lambda-\mu},
\end{equation}
where $m_{\lambda-\mu}
=\min_{1\leq i\leq\ell}\langle\lambda-\mu,\omega_i\rangle$. 
\end{definition}

The following theorem is the main result of this paper.\nopagebreak
\begin{theorem}\label{thm:main1}
Let $G$ be an adjoint simple algebraic group. For $\lambda \in X_*(T)^+$, if there exists $\mu\in X_*(T)^+$ such that  $(\lambda, \mu)$ forms a $\theta$-minimal pair, then we have
\[
T_{e_\mu}\Gr_{\preceq\lambda}\subsetneq T_{e_\mu}\Gr^{\rm FM}_{\preceq\lambda}.
\]
In particular, the FM Schubert scheme $\Gr^{\rm FM}_{\preceq\lambda}$ is nonreduced.
\end{theorem}

\begin{proof}
Set $q=\langle\lambda-\mu,\theta\rangle$. By $\theta$-minimality,
$1\leq q\leq m_{\lambda-\mu}$.
Since  $\Gr_{\preceq\lambda} \hookrightarrow \Gr^{\rm FM}_{\preceq\lambda}$ is a closed embedding, we have $T_{e_\mu}\Gr_{\preceq\lambda}\subseteq T_{e_\mu}\Gr^{\rm FM}_{\preceq\lambda}$. By Proposition~\ref{prop:fmtest}, $q\le m_{\lambda-\mu}$ gives 
\begin{equation}\label{eq:incl}
 \xi_{\mu,q}\in T_{e_\mu}\Gr^{\rm FM}_{\preceq\lambda}.\end{equation}

\smallskip\noindent
We first consider the simply-laced case.  In view of Corollary~\ref{cor:polo}, $1\leq\langle\lambda-\mu,\theta\rangle\leq q$ gives $\xi_{\mu,q}\notin T_{e_\mu}\Gr_{\preceq\lambda}$.

\smallskip\noindent
We now assume $G$ is a non simply-laced adjoint simple group.  We can realize $G=\tilde{G}^\sigma$, for some simply-laced and adjoint simple group  $\tilde{G}$, and  a pinned diagram automorphism $\sigma$, using the following fixed point realization:
\[
D_{n+1}\longrightarrow B_n\ (n\geq3),\qquad
E_6\longrightarrow F_4,\qquad
D_4\longrightarrow G_2,\qquad
A_{2n-1}\longrightarrow C_n\ (n\geq2),
\]
cf.\,\cite[Section~2]{HongMV}.
Here $B_2$ is included as $C_2$ with the nodes interchanged;
the folding of type $A_{2n}$ is excluded. Note that for these $\tilde{G}$ considered,  $\tilde{G}^\sigma$ is connected and adjoint, cf.\,\cite[Propositions~1.30(1)]{DM}.

Let $\tilde{T}\subset\tilde{G}$ be
the pinned maximal torus and put $T=\tilde{T}^\sigma$.
Then
$X_*(T)=X_*(\tilde{T})^\sigma$, and for an orbit $\imath$ of vertices
of the Dynkin diagram of $\tilde{G}$,
\begin{equation}\label{eq:fold-coroot}
\alpha_{\imath}^\vee
=\sum_{j\in\imath}\tilde{\alpha}_j^\vee;
\end{equation}
see \cite[Proposition~2.1]{HongMV} or \cite[Lemma~2.2]{HongShen}. The coweights
 $\lambda$ and $\mu$ can be regarded as $\sigma$-invariant dominant coweights of $\tilde{G}$, and
their dominance order agrees in the two groups. Write $\lambda-\mu=\sum \tilde{a}_i \tilde{\alpha}_i^\vee$.  Let $\tilde{m}_{\lambda-\mu}$ denote the minimal coefficient among $\tilde{a}_i$. 
In view of \eqref{eq:fold-coroot}, $m_{\lambda-\mu}=\tilde{m}_{\lambda-\mu}$.

Let $\tilde{\theta}$ be the highest root of $\tilde{G}$, and put
$\tilde{\g}=\operatorname{Lie}(\tilde{G})$. For these foldings,
$\tilde{\theta}|_T=\theta$, and $\sigma$ acts trivially on
$\tilde{\g}_{\tilde{\theta}}$; see \cite[Lemma~4.3]{HongMV}.
Choose $x_\theta=x_{\tilde{\theta}}$, and let
$\iota:\Gr_G\hookrightarrow\Gr_{\tilde{G}}$ be the natural closed immersion.
By \eqref{eq:infinitesimal}, its differential sends
$\xi_{\mu,q}\in T_{e_\mu}\Gr_G$ to the tangent vector
$\tilde{\xi}_{\mu,q}\in T_{e_\mu}\Gr_{\tilde{G}}$ represented by
$(1+\epsilon x_{\tilde{\theta}}(-q))e_\mu$. Thus
\[
d\iota_{e_\mu}(\xi_{\mu,q})=\tilde{\xi}_{\mu,q}.
\]
Since $\langle\lambda-\mu,\tilde{\theta}\rangle
=\langle\lambda-\mu,\theta\rangle$, applying Corollary~\ref{cor:polo}
for $\tilde{G}$ gives
$\tilde{\xi}_{\mu,q}\notin T_{e_\mu}\Gr^{\tilde{G}}_{\preceq\lambda}$.
On the other hand, the closed immersion $\iota$ induces an isomorphism 
\[
\left (  \bigl(
\Gr_{\preceq\lambda}^{\widetilde G}
\bigr)^\sigma \right )_{\mathrm{red}}  = \Gr_{\preceq\lambda},
\]
where $ \bigl(\Gr_{\preceq\lambda}^{\widetilde G}
\bigr)^\sigma$ denotes the schematic $\sigma$ fixed-point subscheme, 
cf.\,\cite[Corollary 3.3]{HongMV}.
Thus, this induces an embedding of  tangent spaces 
\begin{equation}\label{eq:fold-tangent-inclusion}
T_{e_\mu}\Gr_{\preceq\lambda}=T_{e_\mu} \left (  \bigl(
\Gr_{\preceq\lambda}^{\widetilde G}
\bigr)^\sigma \right )_{\mathrm{red}} 
\hookrightarrow T_{e_\mu} \left ( \bigl(  
\Gr_{\preceq\lambda}^{\widetilde G}
   \bigr)  ^\sigma    \right)  \simeq
\bigl(T_{e_\mu}\Gr^{\tilde{G}}_{\preceq\lambda}\bigr)^\sigma.
\end{equation}
Since $\tilde{\xi}_{\mu,q}\notin \bigl(T_{e_\mu}\Gr^{\tilde{G}}_{\preceq\lambda}\bigr)^\sigma$, it follows that $\xi_{\mu,q}\notin T_{e_\mu}\Gr_{\preceq\lambda}$.  In view of Lemma \ref{lem:FM_red} and \eqref{eq:incl}, we have $T_{e_\mu}\Gr_{\preceq\lambda}\varsubsetneq T_{e_\mu}\Gr^{\rm FM}_{\preceq\lambda}$. It follows that $\Gr^{\rm FM}_{\preceq\lambda}$ is not reduced.
\end{proof}

\section{$\theta$-minimal pairs and their classification}\label{sect:adm}

In this section, we study the basic properties of $\theta$-minimal pairs and
classify the dominant coweights which admit such pairs in every simple type.
We also construct explicit families in types $B,D,E,F,G$.
In type $D_n$, we show that no such pair has $\mu=0$, and give examples of
nonreduced FM Schubert schemes whose tangent spaces at the base point agree
with those of the corresponding Schubert varieties.

We keep the root and coweight notation of Section~\ref{sect:main_thm}
and assume throughout this section that $G$ is adjoint.
\subsection{Properties of $\theta$-minimal pairs and the classification theorem}
\label{subsec:non-simply-laced-examples}

Table~\ref{table:highest-root} gives the highest root in each type.
We use Bourbaki numbering throughout; in particular,
$\alpha_1,\alpha_2$ are long in type $F_4$, and $\alpha_1$ is short
in type $G_2$.
\medskip
\begin{table}[ht]
\centering
\renewcommand{\arraystretch}{1.25}
\setlength{\tabcolsep}{4pt}
\footnotesize
\begin{tabularx}{\linewidth}{@{}
    >{\centering\arraybackslash}p{0.12\linewidth}
    >{\centering\arraybackslash}p{0.22\linewidth}
    >{\raggedright\arraybackslash}X
    @{}}
\toprule
Type
& Fundamental weight expression of $\theta$
& Simple root expression of $\theta$ \\
\midrule

$A_r$, $r\ge 1$
& $\omega_1+\omega_r$
& $\sum_{i=1}^r\alpha_i$ \\

$B_2$
& $2\omega_2$
& $\alpha_1+2\alpha_2$ \\

$B_n$, $n\ge3$
& $\omega_2$
& $\alpha_1+2\alpha_2+\cdots+2\alpha_n$ \\

$C_n$, $n\ge2$
& $2\omega_1$
& $2\alpha_1+\cdots+2\alpha_{n-1}+\alpha_n$ \\

$D_n$, $n\ge 4$
& $\omega_2$
& $\alpha_1+2\alpha_2+\cdots+2\alpha_{n-2}
   +\alpha_{n-1}+\alpha_n$ \\

$E_6$
& $\omega_2$
& $\alpha_1+2\alpha_2+2\alpha_3+3\alpha_4
   +2\alpha_5+\alpha_6$ \\

$E_7$
& $\omega_1$
& $\begin{aligned}
   &2\alpha_1+2\alpha_2+3\alpha_3+4\alpha_4
   +3\alpha_5+2\alpha_6+\alpha_7
   \end{aligned}$ \\

$E_8$
& $\omega_8$
& $\begin{aligned}
   &2\alpha_1+3\alpha_2+4\alpha_3+6\alpha_4
   +5\alpha_5+4\alpha_6+3\alpha_7+2\alpha_8
   \end{aligned}$ \\
$F_4$
& $\omega_1$
& $2\alpha_1+3\alpha_2+4\alpha_3+2\alpha_4$ \\

$G_2$
& $\omega_2$
& $3\alpha_1+2\alpha_2$ \\
\bottomrule
\end{tabularx}
\caption{The highest root $\theta$ for the irreducible root systems.}
\label{table:highest-root}
\end{table}

For the types in which $\theta$ is a fundamental weight, define
\begin{equation}\label{eq:node}
s=\begin{cases}
2,&\g=B_n\ (n\ge3),\ D_n,\ E_6,\text{ or }G_2,\\
1,&\g=F_4\text{ or }E_7,\\
8,&\g=E_8.
\end{cases}
\end{equation}
Thus $\theta=\omega_s$ in these types. Write
$\lambda-\mu=\sum_{i=1}^\ell a_i\alpha_i^\vee$, where
$a_i=\langle\lambda-\mu,\omega_i\rangle$.

\begin{lemma}\label{lem:criterion}
There are no $\theta$-minimal pairs in types $A$ and $C$, or in type
$B_2$. In every remaining type, a pair $\mu\prec\lambda$ in
$X_*(T)^+$ is $\theta$-minimal if and only if
\begin{equation}\label{eq:criterion}
a_s=\min_i a_i\geq1.
\end{equation}
\end{lemma}
\begin{proof}
Put $k:=\langle\lambda-\mu,\theta\rangle$. For $A_\ell$, $\ell\ge2$, $k=a_1+a_\ell$. If $(\lambda, \mu)$ is $\theta$-minimal we would have inequalities $1\le a_1+a_\ell\le\min_ia_i\le a_1$, so $a_\ell\le0$, contradicting $a_\ell\ge\min_ia_i\ge k\ge1$. For $A_1$, $k=2a_1\le a_1$ forces $a_1\le0$, contradicting $k\ge1$. 

For $C_n$, the same argument applies since $k=2a_1$; for $B_2$, it applies with $a_2$ in place of $a_1$.

For the remaining types, observe from the above table,  \eqref{eq:question}  is equivalent to the following condition:
$1\leq  a_s\le\min_i a_i$.
\end{proof}

The following lemma shows that the $\theta$-minimal pairs are closed under translation, addition and scaling.
\begin{lemma}\label{lem:operations}
If $(\lambda,\mu)$ is $\theta$-minimal and $\nu\in X_*(T)^+$, then
$ (\lambda+\nu,\mu+\nu)$
is $\theta$-minimal. If $(\lambda_1,\mu_1)$ is another $\theta$-minimal pair
for the same root system, then
$ (\lambda+\lambda_1,\mu+\mu_1)$
is $\theta$-minimal. In particular, $(k\lambda,k\mu)$ is $\theta$-minimal for
every integer $k\ge1$.
\end{lemma}
\begin{proof}
Any $\theta$-minimal pair $(\lambda,\mu)$ translated by an element in $X_*(T)^+$ is clearly still $\theta$-minimal. 
For addition, it follows from the elementary fact that  $a_i\ge a_s$ and $a_i'\ge a_s'$ implies that
$a_i+a_i'\ge a_s+a_s'$. Therefore, the $
\theta$-minimal pairs are also closed under rescaling. 
\end{proof}

The following theorem gives a complete classification of dominant coweights $\lambda$ that admit a $\theta$-minimal pair. By Theorem~\ref{thm:main1}, their associated FM Schubert schemes $\Gr_{\preceq \lambda}^{\rm{FM}}$ are all nonreduced.
\begin{theorem}\label{thm:classification-all}
Write $\lambda=\sum_i b_i\omega_i^\vee\in X_*(T)^+$. There exists $\mu\in X_*(T)^+$ such that $\mu\prec\lambda$ and $(\lambda,\mu)$ is $\theta$-minimal exactly in the following cases:
\begin{center}
\small
\renewcommand{\arraystretch}{1.18}
\begin{tabularx}{0.98\textwidth}{@{}l@{\hspace{2em}}X@{}}
\toprule
Type & Condition on $\lambda$ \\
\midrule
$A_n$, $C_n$, $B_2$ & none \\
$B_n$, $n\ge3$ & $b_1\ge1$ and $\sum_{i=3}^n b_i\ge1$ \\
$D_n$, $n\ge5$ & $b_1\ge1$ and $2\sum_{i=3}^{n-2}b_i+b_{n-1}+b_n\ge2$ \\
$D_4$ & $b_1,b_3,b_4\ge1$ \\
$E_6$ & $5b_1+4b_3+3b_4+2b_5+b_6\ge3$ and $b_1+2b_3+3b_4+4b_5+5b_6\ge3$ \\
$E_7$ & $b_2+b_7+2(b_3+b_4+b_5+b_6)\ge2$ \\
$E_8$ & $\lambda\neq0$ \\
$F_4$ & $b_2+b_3+b_4\ge1$ \\
$G_2$ & $b_1\ge1$ \\
\bottomrule
\end{tabularx}
\end{center}
\end{theorem}

The types $A$, $C$, and $B_2$ follow from Lemma~\ref{lem:criterion}. The type $D$ cases are proved in Theorem~\ref{prop:minuscule-classification} below. We prove the $E_8$ case here and the remaining cases in Appendix~\ref{app:classification}.

In type $E_8$, the coweight and coroot lattices coincide.  Hence $\omega_i-\theta$ is a nonnegative integral sum of simple roots for every $i$. For any nonzero $\lambda\in X_*(T)^+$, dominance gives
\[
1\le\langle\lambda,\theta\rangle
\le\min_i\langle\lambda,\omega_i\rangle.
\]
Thus, $(\lambda,0)$ is $\theta$-minimal. This proves the $E_8$ case.

\begin{corollary}\label{cor:E8-all}
Assume $G$ is of type $E_8$. Then $\Gr^{\rm FM}_{\preceq\lambda}$ is reduced if and only if $\lambda=0$. In other words, every nontrivial FM Schubert scheme of type $E_8$ is nonreduced.
\end{corollary}
\begin{proof}
By the preceding argument, every nonzero dominant coweight $\lambda$ gives a $\theta$-minimal pair $(\lambda,0)$, so Theorem~\ref{thm:main1} implies nonreducedness. For $\lambda=0$, the FM Schubert scheme is a point.
\end{proof}

\begin{remark}
One can easily check that the examples considered in \cite{BHY} are $\theta$-minimal:
\[  E_6:   (\omega_1^\vee+\omega_6^\vee, 0);   \; E_7:  (\omega_6^\vee, 0);   \; E_8: (\omega_8^\vee, 0) .   \]
\end{remark}

Theorems~\ref{thm:classification-all} and~\ref{thm:main1} also give the following examples in types $B$, $F$, and $G$.
\begin{proposition}\label{prop:fold-examples}
The pair $(\lambda, \mu)$ in each row of the following table is a $\theta$-minimal pair, and thus 
$\Gr^{\mathrm{FM}}_{\preceq\lambda}$ is nonreduced.
\[
\begin{array}{c|c|c|c}
G&\lambda&\mu&a_s\\ \hline
B_n,\ n\geq3&\omega_1^\vee+\omega_3^\vee&\omega_2^\vee&1\\
F_4&\omega_4^\vee&0&2\\
G_2&\omega_1^\vee&\omega_2^\vee&1.
\end{array}
\]
Note that these $\lambda$ and $\mu$ lie in the coroot lattice of $G$.
\end{proposition}

\begin{proof}
The proof of Theorem~\ref{thm:classification-all} in
Appendix \ref{app:classification} gives the displayed $\theta$-minimal pairs and the values of $a_s$.
The nonreducedness follows from Theorem~\ref{thm:main1}.

In type $B_n$, we have
\[
\mu=\omega_2^\vee
=\alpha_1^\vee+2\sum_{i=2}^{n-1}\alpha_i^\vee+\alpha_n^\vee
\in Q^\vee.
\]
Since $\lambda-\mu\in Q^\vee$, both coweights lie in the coroot lattice.
In types $F_4$ and $G_2$, the coweight and coroot lattices coincide.
\end{proof}

Theorem~\ref{thm:classification-all} specializes to the following
classification for sums of minuscule coweights.
\begin{corollary}\label{cor:minuscule-sums}
Let $\lambda\in X_*(T)^+$ be a sum of dominant minuscule coweights.
There exists $\mu\in X_*(T)^+$ such that $\mu\prec\lambda$ and
$(\lambda,\mu)$ is $\theta$-minimal exactly in the following cases:
\begin{center}
\small
\renewcommand{\arraystretch}{1.2}
\begin{tabularx}{0.98\textwidth}{@{}l@{\hspace{1.5em}}l@{\hspace{1.5em}}X@{}}
\toprule
Type & $\lambda$ & Condition \\
\midrule
$D_n$, $n\ge5$ & $r\omega_1^\vee+s\omega_{n-1}^\vee+t\omega_n^\vee$
& $r\ge1$ and $s+t\ge2$ \\
$D_4$ & $r\omega_1^\vee+s\omega_3^\vee+t\omega_4^\vee$
& $r,s,t\ge1$ \\
$E_6$ & $r\omega_1^\vee+s\omega_6^\vee$
& $r+s\ge3$ or $r=s=1$ \\
$E_7$ & $r\omega_7^\vee$ & $r\ge2$ \\
\bottomrule
\end{tabularx}
\end{center}
Here $r,s,t\in\mathbb Z_{\ge0}$. No sum of dominant minuscule coweights
admits a $\theta$-minimal pair in the other types.
In every listed case, $\Gr^{\rm FM}_{\preceq\lambda}$ is nonreduced.
\end{corollary}
\begin{proof}
A fundamental coweight $\omega_i^\vee$ is minuscule if and only if the
coefficient of $\alpha_i$ in $\theta$ is $1$.
Thus Table~\ref{table:highest-root} determines the allowed indices in
$\lambda=\sum_i b_i\omega_i^\vee$. Setting all other $b_i$ equal to zero
in Theorem~\ref{thm:classification-all} gives the stated conditions.
In type $E_6$, the two inequalities become $5r+s\ge3$ and $r+5s\ge3$;
for nonnegative integers $r,s$, these are equivalent to $r+s\ge3$ or
$r=s=1$. The nonreducedness assertions follow from
Theorem~\ref{thm:main1}.
\end{proof}

\subsection{Classification in type $D_n$ }
\label{sec:minuscule-examples}

In the following theorem, we give a classification of the dominant coweights in type
$D_n$ which admit a $\theta$-minimal pair. In particular, for sums of minuscule
coweights the two families considered in \cite[Theorem~4.2.3]{KPZ} form part of the
complementary range.

\begin{samepage}
\begin{theorem}\label{prop:minuscule-classification}
Assume that $G$ is of type $D_n$, with $n\geq 4$, and write
\[
\lambda=\sum_{i=1}^n b_i\omega_i^\vee\in X_*(T)^+.
\]
\begin{enumerate}
\item[(i)] If $n\geq 5$, there exists $\mu\in X_*(T)^+$ such that
$\mu\prec\lambda$ and $(\lambda,\mu)$ is $\theta$-minimal if and only if
\[
b_1\geq 1
\qquad\text{and}\qquad
2\sum_{i=3}^{n-2}b_i+b_{n-1}+b_n\geq 2.
\]

\item[(ii)] If $n=4$, such a $\mu$ exists if and only if 
$b_1,b_3,b_4\geq 1$.
\end{enumerate}
\end{theorem}
\end{samepage}

\begin{proof}
First we prove $(\mathrm{i})$, so assume $n \geq 5$. Recall that $\theta = \omega_2$. By Lemma \ref{lem:operations}, to prove sufficiency, it is enough to find a $\theta$-minimal pair $(\lambda_0, \mu_0)$ with $\lambda- \lambda_0 \in X_*(T)^+$. Under the assumptions in $(\mathrm{i})$, we examine the following four cases.

\begin{enumerate}
\item Case 1: $b_j \geq 1$ for some $3 \leq j \leq n-2$. Take $(\lambda_0,\mu_0)=(\omega_1^\vee+\omega_{j}^\vee, \omega_{j-1}^\vee)$. We have
\[
\omega_1^\vee+\omega_j^\vee-\omega_{j-1}^\vee=\alpha_1^\vee + \dots +\alpha_{j-1}^\vee+2\alpha_j^\vee +\dots +2 \alpha_{n-2}^\vee+\alpha_{n-1}^\vee+\alpha_n^\vee.
\]
\item Case 2: $b_3= \dots =b_{n-2}=0$ and $b_{n-1}, b_{n} \geq 1$. Take $(\lambda_0,\mu_0)=(\omega_1^\vee + \omega_{n-1}^\vee+\omega_n^\vee, \omega_{n-2}^\vee)$. We compute
\[
\omega_1^\vee+\omega_{n-1}^\vee+\omega_{n}^\vee-\omega_{n-2}^\vee= \sum_{i=1}^{n} \alpha_i^\vee.
\]
\item Case 3: $b_3 = \dots =b_{n-2}=0$ and $b_{n-1} \geq 2, b_{n}=0$. Take $(\lambda_0,\mu_0)=(\omega_1^\vee+2 \omega_{n-1}^\vee, \omega_{n-3}^\vee)$. We compute
\[
\omega_1^\vee+2\omega_{n-1}^\vee- \omega_{n-3}^\vee=\alpha_1^\vee + \dots +\alpha_{n-3}^\vee+2 \alpha_{n-2}^\vee+2\alpha_{n-1}^\vee+\alpha_n^\vee.
\]
\item Case 4: $b_3= \dots =b_{n-2}=0$ and $b_{n-1}=0, b_n \geq 2$. This case is almost the same as Case 3, but with $(\lambda_0, \mu_0)= (\omega_1^\vee+2 \omega_n^\vee, \omega_{n-3}^\vee)$.
\end{enumerate}
In each case, $\lambda-\lambda_0\in X_*(T)^+$, and the simple-coroot coefficients of $\lambda_0-\mu_0$ satisfy $a_2=\min_i a_i=1$. Thus Lemma \ref{lem:operations} gives the required pair.

We now prove the necessity. Assume that $(\lambda,\mu)$ is a $\theta$-minimal pair and write $\lambda-\mu=\sum_{i=1}^{n}a_i\alpha_i^\vee$. By Lemma \ref{lem:criterion}, $a_i\ge a_2\ge1$ for every $i$. Since $\mu$ is dominant,
\[
0\le b_1-2a_1+a_2\le b_1-a_2,
\]
so $b_1\ge1$.

Assume now for the sake of contradiction that
\[
2\sum_{i=3}^{n-2}b_i+b_{n-1}+b_n\le1.
\]
Then $b_3=\cdots=b_{n-2}=0$ and $b_{n-1}+b_n\le1$. Dominance at nodes $n-1,n$ gives
\[
2(a_{n-1}+a_n-a_{n-2})\le b_{n-1}+b_n\le1.
\]
Since the $a_i$ are integers, $a_{n-1}+a_n\le a_{n-2}$, and hence $a_{n-2}\ge2a_2$. On the other hand, successively applying dominance at nodes $n-2,n-3,\ldots,3$ gives
\[
a_{n-2}\le a_{n-3}\le\cdots\le a_2,
\]
a contradiction. This proves $(\mathrm{i})$.

We turn to proving part $(\mathrm{ii})$, so we are working in type $D_4$. Assume $(\lambda,\mu)$ is $\theta$-minimal, $\lambda= \sum b_i \omega_i^\vee$ and $\lambda-\mu= \sum a_i \alpha_i^\vee$. The $\theta$-minimality condition implies $a_2=\min_i a_i\ge1$, and dominance of $\mu$ gives
\[
0\le b_i-2a_i+a_2\le b_i-a_2, \qquad i\in\{1,3,4\},
\]
so we must have $b_1,b_3,b_4\ge1$. Conversely, $(\omega_1^\vee+\omega_3^\vee+\omega_4^\vee, \theta^\vee)$ is $\theta$-minimal, and Lemma \ref{lem:operations} shows that any $\lambda= b_1\omega_1^\vee+b_2\omega_2^\vee+b_3 \omega_3^\vee+b_4\omega_4^\vee$ with $b_1, b_3, b_4 \geq 1$ admits a $\theta$-minimal pair $(\lambda, \mu)$.
\end{proof}

The following family shows that, in applying the
$\theta$-minimality criterion in type $D_n$, one cannot in general
restrict to $\mu=\theta^\vee$. Indeed, there are $\theta$-minimal pairs
$(\lambda,\mu)$ with $\mu\succ\theta^\vee$ for which
$(\lambda,\theta^\vee)$ is not $\theta$-minimal.

\begin{proposition}\label{prop:above-theta}
Assume that $G$ is of type $D_n$, with $n\geq4$, and set
\[
\lambda_0=
\begin{cases}
\omega_1^\vee+\omega_3^\vee+\omega_4^\vee, & n=4,\\[2mm]
\omega_1^\vee+\omega_3^\vee, & n\geq5.
\end{cases}
\]
For every integer $k\geq2$, define
\[
\mu=k\theta^\vee,
\qquad
\lambda=\lambda_0+(k-1)\theta^\vee.
\]
Then,  $\theta^\vee\prec\mu\prec\lambda$, and the pair $(\lambda,\mu)$ is $\theta$-minimal, whereas
$(\lambda,\theta^\vee)$ is not $\theta$-minimal.
\end{proposition}

\begingroup
\postdisplaypenalty=10000
\begin{proof}
By the proof of Theorem~\ref{prop:minuscule-classification}, $(\lambda_0,\theta^\vee)$ is $\theta$-minimal.
Since
\[
(\lambda,\mu)
=
\bigl(\lambda_0+(k-1)\theta^\vee,\,
      \theta^\vee+(k-1)\theta^\vee\bigr),
\]
Lemma~\ref{lem:operations} implies that $(\lambda,\mu)$ is
$\theta$-minimal. 

It remains to show that $(\lambda,\theta^\vee)$ is not $\theta$-minimal.
For every $n\geq4$, we have
\[
\lambda_0=2\theta^\vee-\alpha_2^\vee,
\qquad
\lambda-\theta^\vee=k\theta^\vee-\alpha_2^\vee.
\]
Writing $\lambda-\theta^\vee=\sum_i a_i\alpha_i^\vee$, we obtain
\[
a_2=2k-1>k=a_1
\qquad(k\geq2).
\]
By Lemma~\ref{lem:criterion}, $(\lambda,\theta^\vee)$ is not
$\theta$-minimal.
\end{proof}
\endgroup

\subsection{Tangent spaces at the base point in type $D$}
\label{subsect:exmaple_theta}

\begin{proposition}\label{prop:no-basepoint-D}
Assume that $G$ is of type $D_n$, with $n\geq 4$.
There is no $\theta$-minimal pair of the form $(\lambda,0)$.
Equivalently, every $\theta$-minimal pair $(\lambda,\mu)$ in type $D_n$
satisfies $\mu\neq 0$.
\end{proposition}

\begin{proof}
Suppose that $(\lambda,0)$ is $\theta$-minimal. We have
\[
\theta-\omega_1
=\alpha_2+\cdots+\alpha_{n-2}
+\frac12(\alpha_{n-1}+\alpha_n).
\]
On the other hand, $\lambda$ is dominant, so
$\langle\lambda,\alpha_i\rangle\geq0$ for all $i$. Thus $\theta$-minimality gives
\[
\langle\lambda,\theta\rangle
\leq\langle\lambda,\omega_1\rangle
\leq\langle\lambda,\theta\rangle.
\]
It follows that $\langle\lambda,\alpha_i\rangle=0$ for every $i\geq2$.
Hence $\lambda=b_1\omega_1^\vee$ for some $b_1>0$. But
\[
\langle\lambda,\omega_{n-1}\rangle
=\frac{b_1}{2}<b_1=\langle\lambda,\theta\rangle,
\]
contradicting $\theta$-minimality. Hence, there is no $\theta$-minimal pair of the form
$(\lambda,0)$.
\end{proof}

Thus, our $\theta$-minimality criterion cannot detect the nonreducedness of $\mathrm{Gr}_{\preceq \lambda}^{\rm{FM}}$ at the base point $e$ in type $D_n$. 


We now collect some examples in type $D_n$ where $(\lambda, \theta^\vee)$ is $\theta$-minimal, but $T_e \Gr_{\preceq \lambda}=T_e \Gr_{\preceq \lambda}^{\rm{FM}}$. Using \cite[Table 1, p.69]{Hu}, one easily checks that the following pairs $(\lambda,\mu)$ are $\theta$-minimal.
\begin{table}[H]
\centering
\begin{tabular}{c|c|c|c}
\text{Type} & $\lambda $& $\mu$ & $a_s$ \\
\hline
$D_4$ & $\omega_1^\vee+\omega_3^\vee+\omega_4^\vee$ & $\theta^\vee$ & $1$ \\
$D_n,\ n\geq5 $& $\omega_1^\vee+\omega_3^\vee$ & $\theta^\vee$ & 1 \\
\end{tabular}
\caption{Pairs for Theorem~\ref{lem:nonbase}.}
\label{thm4.5table}
\end{table}
Note that the coweights in Table \ref{thm4.5table} are in the coroot lattice, so $e \in \Gr_{\preceq \lambda}$. From the following theorem, we see it is necessary to examine the tangent spaces at non-base points.
\begin{theorem}\label{lem:nonbase}
Let $\lambda$ be a coweight as in Table \ref{thm4.5table} above.
Despite $(\lambda,\theta^\vee)$ being $\theta$-minimal, we have equality of tangent spaces at the base point:
$$T_e \Gr_{\preceq \lambda}= T_e\Gr^{\rm FM}_{\preceq \lambda}= \g t^{-1}\oplus \g t^{-2}.$$
\end{theorem}
\begin{proof}
We first note that
\[
\lambda= \begin{cases}
2\alpha_1^\vee+3\alpha_2^\vee+2\alpha_3^\vee+2\alpha_4^\vee,  \quad \g=D_4\\
2\alpha_1^\vee+3\alpha_2^\vee +4\sum_{i=3}^{n-2}\alpha_i^\vee +2\alpha_{n-1}^\vee+2\alpha_n^\vee, \quad \g=D_n, \,n\geq 5
\end{cases}
\]
By \cite[Theorem 4.5]{BHY},  $  T_{e} \Gr^{\rm FM}_{\preceq \lambda} = \mathfrak{g}t^{-1} \oplus \mathfrak{g}t^{-2}$. Thus,  it suffices to show that $T_e\mathrm{Gr}_{\preceq \lambda}= \mathfrak{g}t^{-1} \oplus \mathfrak{g}t^{-2}$. First note that $\Gr_{\preceq \theta^\vee}\subset \Gr_{\preceq \lambda}$. By \cite[Proposition 3.2]{BHY}, $\g t^{-1}= T_e\Gr_{\preceq \theta^\vee}\subset  T_e\Gr_{\preceq \lambda}$. What remains is to show that $\mathfrak{g}t^{-2} \subset T_e \Gr_{\preceq \lambda}$ as well.

For the following Frenkel-Kac calculation only, we identify weights and
coweights via the normalized invariant form and use the notation of
Section~\ref{sec:FK}.
We will show that $x_{\theta}(-2)v_0 \in D(1,\lambda)$; by Polo's criterion, this will imply $x_{\theta}(-2) \in T_e \Gr_{\preceq \lambda}$. As a preliminary step we will show certain extremal vectors belong to $D(1,\lambda)$.

Set
$\beta=\theta-\alpha_2$. Then,
both $\beta$ and $\alpha_2$ are roots, and we have the following identities:
\[\beta+\alpha_2=\theta, \text{ as well as } \lambda=\theta+\beta,\]
where we identify $\lambda$ with an element in the weight lattice.

We also have the following pairings:
\[
\qquad (\theta\mid\beta)=(\theta\mid\alpha_2)=1.
\qquad
\]

Since $\lambda=2\theta-\alpha_2$ and
$(\lambda\mid\theta)=3$, we have  $s_\theta(-\lambda)=\theta+\alpha_2$.  Note that we always have $w_0(\theta)=-\theta$.  When $\g=D_4$,  $w_0=-1$; when $\g=D_n, n\geq 5$,  the involution $-w_0$ fixes node 2 in the Dynkin diagram. Thus, for any $D_n$, we have $w_0(\alpha_2)=-\alpha_2$.  Therefore,  $w_0\lambda=-\lambda$.  Note that $D(1,\lambda)$  is $\mathfrak{g}$-stable, hence it is stable under representatives of $W$ in $G$. Since $D(1,\lambda)$ contains the extremal vector $v_{-\lambda}$, it follows that the following extremal vectors also belong to $D(1, \lambda)$:
\begin{equation}\label{eq:v_D}
v_{\lambda}=v_{\theta+\beta},\,
v_{\theta+\alpha_2}\in D(1,\lambda).
\end{equation}
For $\gamma= \beta$ or $\gamma= \alpha_2$, we can compute $x_{-\gamma}(1)v_{\theta+\gamma}$ explicitly using Formula \eqref{eq:VOA-extrem-vec} and we obtain
\[
x_{-\gamma}(1)v_{\theta+\gamma}
=-\epsilon(-\gamma,\theta+\gamma)\,\gamma(-1)v_\theta.
\]
Since $v_{\theta+\beta}$ and $v_{\theta+\alpha_2}$ are both in $D(1, \lambda)$ and since $x_{-\gamma}(1)\in\g[t]$, both
$\beta(-1)v_\theta$ and $\alpha_2(-1)v_\theta$ lie in $D(1,\lambda)$. Therefore
\[
\theta(-1)v_\theta
=\beta(-1)v_\theta+\alpha_2(-1)v_\theta\in D(1,\lambda).
\]
On the other hand, using Formula \eqref{eq:VOA-extrem-vec}, we have  $\theta(-1)v_{\theta}=x_{\theta}(-2)v_0$.
By Polo's criterion (Theorem~\ref{thm:polo}), $x_\theta(-2)\in T_e\Gr_{\preceq\lambda}$.

Since $T_{e}\Gr_{\preceq \lambda}$ is $\g$-stable and since $x_{\theta}(-2) \in T_e \Gr_{\preceq \lambda}$, we conclude that
$\g t^{-2}\subset T_e \Gr_{\preceq \lambda}$, and
thus $T_e \Gr_{\preceq \lambda}= T_e\Gr^{\rm FM}_{\preceq \lambda}= \g t^{-1}\oplus \g t^{-2}$.
\end{proof}
\begin{remark}
This theorem answers a question raised by Alex Weekes: does the equality
$T_e\Gr_{\preceq\lambda}= T_e\Gr^{\mathrm{FM}}_{\preceq\lambda}$
imply that the FM Schubert scheme
$\Gr^{\mathrm{FM}}_{\preceq\lambda}$ is reduced?
The answer is no for simply-connected groups $G$. It is worth mentioning that, by \cite[Proposition~5.1]{KMW}, the reducedness of
$\Gr^{\mathrm{FM}}_{\preceq\lambda}$ in a neighborhood of the base point $e$ does imply its global reducedness.
\end{remark}

\appendix
\addtocontents{toc}{\protect\setcounter{tocdepth}{1}}
\section{The remaining classification cases}\label{app:classification}
In this Appendix, we prove the remaining  cases of Theorem~\ref{thm:classification-all}. By Lemma \ref{lem:criterion}, there is no $\theta$-minimal pair for type $A,C$ and $B_2$. The type $D_n$ case has been proved in Theorem~\ref{prop:minuscule-classification}, and type $E_8$ is demonstrated right after Theorem~\ref{thm:classification-all}. 

In the following case by case arguments, we use Lemma~\ref{lem:operations} freely: once a $\theta$-minimal pair is found for a dominant coweight $\lambda_0$, every $\lambda$ with $\lambda-\lambda_0\in X_*(T)^+$ also admits one. Recall the setup in Theorem~\ref{thm:classification-all},  we write $\lambda = \sum b_i \omega_i^{\vee}$. 

\subsection{Type $B_n$}\label{app:classification-B}
For type $B_n$,  we have $\theta=\omega_2$. We first prove the sufficiency of the criterion. For $3\le j\le n$, the pair
\[
(\lambda_0,\mu_0)=(\omega_1^\vee+\omega_j^\vee,\omega_{j-1}^\vee)
\]
is $\theta$-minimal, since
\[
\lambda_0-\mu_0
=\alpha_1^\vee+\cdots+\alpha_{j-1}^\vee
+2\alpha_j^\vee+\cdots+2\alpha_{n-1}^\vee+\alpha_n^\vee,
\]
with the middle sum empty when $j=n$. Thus $b_1\ge1$ and $\sum_{i=3}^n b_i\ge1$ is sufficient.

Conversely, assume that $(\lambda, \mu)$ is a $\theta$-minimal pair in type $B_n$, $n \geq 3$, and write $\lambda-\mu=\sum a_i\alpha_i^\vee$. Lemma \ref{lem:criterion} asserts that we must have $a_2=\min_i a_i\ge1$. Dominance of $\mu$ at node $1$ gives
\[
0\le b_1-2a_1+a_2,
\]
so $b_1\ge1$. If $b_3=\cdots=b_n=0$, multiplying the dominance inequality at node $i$ by $i-2$ and summing over $3\le i\le n$ gives $2a_n\le a_2$, contradicting $a_n\ge a_2\ge1$. This proves the type $B_n$ criterion.

\subsection{Type $E_6$}
For $E_6$, we have $\theta=\omega_2$.  Assume $(\lambda, \mu)$ is $\theta$-minimal and write $\lambda-\mu=\sum_i a_i\alpha_i^\vee$. By Lemma \ref{lem:criterion}  we know that $a_2=\min_i a_i\ge1$. Consider the following relation:
\[
3(2\omega_1-\omega_2)
=5\alpha_1+4\alpha_3+3\alpha_4+2\alpha_5+\alpha_6.
\]
Since $\mu$ is a dominant coweight, we must have  $\langle \mu, 3(2 \omega_1-\omega_2)\rangle\geq 0$. It follows that
\[
5b_1+4b_3+3b_4+2b_5+b_6\ge6a_1-3a_2\ge3,
\]
where the second inequality is due to $a_2 = \min_i a_i \geq 1$.
By symmetry, we also have
\[
b_1+2b_3+3b_4+4b_5+5b_6\ge6a_6-3a_2\ge3,
\]
so we see these two inequalities are necessary for $(\lambda, \mu)$ to be $\theta$-minimal.

We claim these inequalities are also sufficient. Their nonnegative integral solutions form the upper ideal generated by the coweights in the first column of the following table and their images under the diagram involution. The last column gives the coefficients of $\lambda_0-\mu_0$ in the simple-coroot basis:
\[
\begin{array}{c|c|c}
\lambda_0&\mu_0&(a_1,\ldots,a_6)\\ \hline
\omega_4^\vee&\omega_2^\vee&(1,1,2,3,2,1)\\
2\omega_3^\vee&\omega_1^\vee+\omega_2^\vee&(1,1,3,3,2,1)\\
3\omega_1^\vee&\omega_2^\vee&(3,1,3,3,2,1)\\
\omega_1^\vee+\omega_3^\vee&\omega_2^\vee&(2,1,3,3,2,1)\\
\omega_1^\vee+\omega_5^\vee&\omega_3^\vee&(1,1,1,2,2,1)\\
\omega_1^\vee+\omega_6^\vee&0&(2,2,3,4,3,2)\\
\omega_3^\vee+\omega_5^\vee&\omega_4^\vee&(1,1,2,2,2,1)
\end{array}
\]
In every row $a_2=\min_i a_i\ge1$, proving the $E_6$ case.

\subsection{Type $E_7$} We first show that if $(\lambda, \mu)$ is a $\theta$-minimal pair then it is necessary that $b_2+b_7+2(b_3+b_4+b_5+b_6) \geq 2$. Assume for the sake of contradiction that the inequality does not hold, so we have $b_3=\cdots=b_6=0$ and $b_2+b_7\le1$. Write $\lambda-\mu = \sum a_i \alpha_i^{\vee}$. If $(\lambda, \mu)$ is $\theta$-minimal, then by Lemma \ref{lem:criterion} we have $1 \leq \min_i a_i = a_1$.
We also express $\mu$ as a sum of fundamental coweights; 
\[\sum b_i \omega_i^{\vee}-\sum a_i \alpha_i^{\vee}=\mu=\sum_i c_i\omega_i^\vee\] with $c_i\in\mathbb Z_{\ge0}.$

We compute $\langle \mu, \alpha_2+\alpha_3+2(\alpha_4+\alpha_5+\alpha_6+\alpha_7) \rangle$ in two different ways, and obtain the following equality relating $a_i$'s, $b_i$'s and $c_i$'s:
 
\begin{equation}
a_1+b_2+2b_7-2a_7=c_2+c_3+2(c_4+c_5+c_6+c_7).
\label{eq:e7}
\end{equation}
Note that by the $\theta$-minimality of $(\lambda, \mu)$, we have $a_7\ge a_1\geq 1$ and then 
$$a_1-2a_7\leq -a_7 \leq -1.$$ Thus, we must have $b_2+b_7=1$, else the left hand side of \eqref{eq:e7} would be negative, whereas the right hand side is manifestly nonnegative. 
 A quick check by cases shows
$0\le c_2+c_3+2(c_4+c_5+c_6+c_7) \le1$. Thus $c_4=\cdots=c_7=0$. Using the vanishing of these coefficients, we compute the pairing \[0=\langle \mu, \alpha_5+2\alpha_6+3\alpha_7 \rangle = a_4-4a_7+3b_7,\]
or in other words $a_4=4a_7-3b_7$. Consequently,
\begin{equation}
	c_2=\langle \mu, \alpha_2 \rangle = b_2-2a_2+a_4=b_2-2a_2-3b_7+4a_7.
	\label{eq:e7-2}
\end{equation}

In either of the two cases ($b_2=0, b_7=1$ or $b_2=1, b_7=0$), the integer on the right is odd. Since we have $0 \leq c_2 + c_3 \leq 1$ and $c_2$ is odd, we see $c_2=1$ and $c_3=0$.

Plugging this back to equation \ref{eq:e7}, we have 
\[a_1+b_2+2b_7=1+2a_7.\]
If $b_2=1$ and $b_7=0$, then $a_1+1=1+2a_7$ has no positive integer solutions $1 \leq a_1 \leq a_7$. The other case $b_2=0, b_7=1$ gives us the equation $a_1+2=1+2a_7$; the only solution when $1 \leq a_1 \leq a_7$ is $a_1=a_7=1$. But then plugging this into \ref{eq:e7-2}, we obtain 
\[1= c_2 = b_2-2a_2-3b_7+4a_7=0-2a_2-3+4.\]
This forces $a_2=0$, which contradicts the $\theta$-minimality condition $1 \leq a_1 \leq a_2$.
Thus the stated inequality is necessary.

We now show the sufficiency. The minimal dominant coweights $\lambda_0=\sum b_i \omega_i^{\vee}$ such that $b_2+b_7+2(b_3+b_4+b_5+b_6) \geq 2$ in $E_7$ are given in the following table, which gives a $\theta$-minimal pair for each:
\[
\begin{array}{c|c|c}
\lambda_0&\mu_0&(a_1,\ldots,a_7)\\ \hline
\omega_3^\vee&\omega_1^\vee&(1,2,3,4,3,2,1)\\
\omega_4^\vee&\omega_3^\vee&(1,2,2,4,3,2,1)\\
\omega_5^\vee&\omega_2^\vee&(1,1,2,3,3,2,1)\\
\omega_6^\vee&0&(2,3,4,6,5,4,2)\\
2\omega_2^\vee&\omega_3^\vee&(1,3,2,4,3,2,1)\\
\omega_2^\vee+\omega_7^\vee&\omega_1^\vee&(1,3,3,5,4,3,2)\\
2\omega_7^\vee&0&(2,3,4,6,5,4,3)
\end{array}
\]
with $a_1=\min_i a_i$ in every row. This proves sufficiency in the $E_7$ case.

\subsection{Type $F_4$}\label{app:classification-F} In type $F_4$, for $\lambda = b_1\omega_1^{\vee}+b_2\omega_2^{\vee}+b_3\omega_3^{\vee}+b_4 \omega_4^{\vee}$, there is a $\theta$-minimal pair $(\lambda,\mu)$ if and only if $b_2+b_3+b_4 \geq 1$.
For $F_4$, we have $\theta=\omega_1$. The pairs
\[
(\omega_2^\vee,\omega_1^\vee),\qquad
(\omega_3^\vee,\omega_2^\vee),\qquad
(\omega_4^\vee,0)
\]
are $\theta$-minimal, since their differences have simple-coroot coefficient vectors
\[
(1,3,2,1),\qquad (1,2,2,1),\qquad (2,4,3,2).
\]
Thus $b_2+b_3+b_4\ge1$ is sufficient. 

To see necessity, assume $\lambda=b_1\omega_1^\vee$ belongs to a $\theta$-minimal pair $(\lambda, \mu)$. Write $\lambda-\mu = \sum a_i \alpha_i^{\vee}$.  Then $a_1=\min_i a_i\ge1$. To obtain a contradiction, we consider the  relation 
$2\omega_4-\omega_1=\alpha_2+2\alpha_3+2\alpha_4$ and pair it against $\lambda-\mu $. 
By the dominance of  $\mu$, we get 
\[
2a_4-a_1
=\langle\lambda-\mu,2\omega_4-\omega_1\rangle = \langle -\mu,2\omega_4-\omega_1\rangle 
\le0,
\]
contradicting the $\theta$-minimality assumption $a_4\ge a_1\geq 1$. Hence the $F_4$ criterion follows.
\subsection{Type $G_2$} For $\lambda = b_1 \omega_1^{\vee} +b_2 \omega_2^{\vee}$, we show that there exists a $\theta$-minimal pair $(\lambda, \mu)$ if and only if $b_1 \geq 1$.
For $G_2$, we have $\theta=\omega_2$ and
\[
\omega_1^\vee-\omega_2^\vee=\alpha_1^\vee+\alpha_2^\vee.
\]
So, $(\omega_1^{\vee}, \omega_2^{\vee})$ is $\theta$-minimal and $b_1\ge1$ is sufficient. 

For necessity, assume $(\lambda, \mu)$ is $\theta$-minimal and $b_1=0$. Dominance of $\mu$ yields $ \langle \mu, \alpha_1 \rangle\geq 0$ and thus $2a_1\le a_2$, whereas $\theta$-minimality gives $a_2\le a_1$ and $a_2\ge1$, a contradiction. This proves the $G_2$ criterion.

\addtocontents{toc}{\protect\setcounter{tocdepth}{2}}
\section*{Disclosure of AI use}
GPT-5.6 Sol and GPT-6 Astra were used as research and writing aids to explore examples, suggest proof strategies for several arguments, and improve exposition.  In particular, AI assistance contributed to finding a counterexample in type $D_4$ corresponding to $(\lambda,\mu)=(\omega_1^\vee+\omega_3^\vee+\omega_4^\vee,\theta^\vee)$. Starting from this example,  the authors subsequently formulated a general  $\theta$-minimality criterion and developed the Frenkel-Kac and Segal-Sugawara approach to the Demazure-module non-inclusion problem.  Any AI-assisted mathematical suggestions were independently verified by the authors. The authors take full responsibility for all mathematical statements and proofs in the manuscript.

\end{document}